# Virtual Transmission Method,

# A New Distributed Algorithm to Solve Sparse Linear System

Fei Wei

Huazhong Yang

Department of Electronic Engineering, Tsinghua University, Beijing, China

Abstract

In this paper, we propose a new parallel algorithm which could work naturally on the parallel computer with arbitrary number of processors. This algorithm is named Virtual Transmission Method (VTM). Its physical background is the lossless transmission line and microwave network. The basic idea of VTM is to insert the virtual transmission lines into the linear system to achieve distributed computing.

VTM is proved to be convergent to solve SPD linear system. Preconditioning method and performance model are presented. Numerical experiments show that VTM is efficient, accurate and stable.

Accompanied with VTM, we bring in a new technique to partition the symmetric linear system, which is named Electric Vertex Splitting (EVS). It is based on Kirchhoff's Current Law from circuit theory. We proved that EVS is feasible to partition any SPD linear system.

**Key words**: Distributed Algorithm, Sparse Linear System, Partitioning, Convergence, Performance Modeling, Preconditioning, Transmission Line, Interconnect, Wire

**New words**: Virtual Transmission Method (VTM), Electric Vertex Splitting (EVS), Transmission Delay Equations (TDE), Virtual Transmission Line (VTL), Interconnect, Wire, Electric Graph, Neighbor-To-Neighbor (N2N), Conformal Splitting Existence

Theory, Impedance Matching, Lines Coupling, Wire Tearing.

# 1. Introduction

The linear system, **Ax** = **b**, is widely encountered in scientific computing. When the coefficient matrix **A** is symmetric-positive-definite (SPD), the linear system is called SPD system, which is extremely common in engineering applications [1, 2]. For example, most of the linear systems generated by the finite element method are SPD systems. Therefore, in many scientific disciplines, solving SPD systems is an inevitable task and the efficiency will be the dominant factor in those fields.

To solve the SPD system, there are two basic approaches, direct methods and iterative methods.

The direct methods are mainly based on the Sparse Cholesky Factorization. In order to efficiently compute the dense submatrices inside the sparse matrix, supernodal method and multifrontal method are used [3].

The representatives of the iterative methods are Conjugate Gradient method (CG) and Multigrid method (MG). CG is based on the Krylov subspace projection. If the preconditioner is properly chosen, the convergence of CG will be fast. MG is efficient for the linear systems generated from the elliptic partial differential equations [4].

All the algorithms mentioned above work well on the traditional single-processor computers, but they would get into trouble on parallel computers [5, 6]. The parallel version of Sparse Cholesky Factorization suffers from the limited concurrency which depends on the distribution of the nonzero elements in the sparse matrix. For the parallel CG, it is difficult to choose a proper preconditioner in a parallel way [4].

Another well known parallel method for large sparse linear system is the Domain Decomposition Method (DDM). DDM refers to a collection of techniques which revolve around the principle of divide and conquer [4]. Schur Complement method, Additive Schwarz method and the Dual-Prime Finite Element Tearing and Interconnection (FETI-DP) method are three commonly-adopted parallel methods of DDM [7].

The Schur Complement method makes use of the master-slave model [8]. This method first partitions the large linear system into a number of subsystems. Then these subsystems are simplified and solved by the slave processors in parallel. After that the simplified results are merged into a new linear system, which is much smaller than the original one. At last this new system is solved by the master processor. This model suffers from the heavy communication overheads imposed on the master processor, especially when the number of slave processors is large. Consequently, the scalability and concurrency of the Schur Complement method is limited.

The Additive Schwarz method is similar to the block Jacobi iteration. For a SPD system, it needs two assumptions to be convergent, and the convergence speed depends on these two assumptions [4].

The FETI-DP method is a scalable method to solve large problems [7, 9]. FETI-DP has to solve a coarse problem. This procedure needs global communication of the residual errors and the concurrency is difficult to explore. Consequently, the parallel efficiency of FETI-DP is affected.

VTM is a new parallel algorithm for large-scale sparse SPD systems. It is inspired by the behavior of transmission lines in the electrical engineering. Although VTM is a distributed iterative algorithm, it is sure to be convergent because of its physical background.

VTM adopts the Neighbor-To-Neighbor (N2N) communication model, which requires only local communication between adjacent processors, as shown in Fig. 1. Because of the N2N model, the communication network of the parallel computer could be simple.

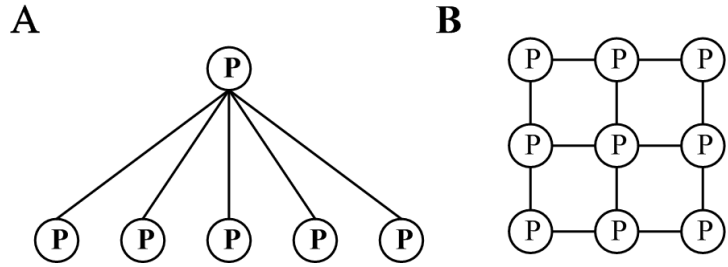

**Figure 1.** Master-slave model Vs. N2N model. (**A**) Master-slave model. (**B**) N2N model.

The paper is organized as follows. Section 2 introduces the basics of transmission line. Section 3 defines the electric graph of the symmetric linear system. Section 4 describes the partitioning technique for the electric graphs. Section 5 details the algorithm of VTM. Section 6 presents the convergence theory for VTM and a basic proof is given in the appendix. Section 7 focuses on the preconditioning of VTM. Section 8 proposes a performance model. Numerical experiments are shown in Section 9. We conclude this work in Section 10.

## 2. Transmission Line

Transmission line is a magic element in electrical engineering. The circuit diagram of a transmission line is illustrated in Fig. 2. The function of the lossless transmission line could be described by the Transmission Delay Equations, as below.

(2.1)
$$\begin{cases} U_1(t) + Z \cdot I_1(t) = U_2(t-\tau) - Z \cdot I_2(t-\tau) \\ U_2(t) + Z \cdot I_2(t) = U_1(t-\tau) - Z \cdot I_1(t-\tau) \end{cases}$$

where $U_1$ and $I_1$ represent the potential and current of Port 1, and $U_2$ and $I_2$ represent those of Port 2. $t$ is the time, and $\tau$ is the propagation delay. $Z$ is the characteristic impedance of the transmission line [18, 19, 20].

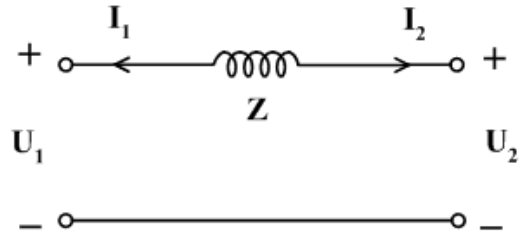

**Figure 2.** The circuit diagram of the transmission line.

Transmission line is always troublesome for integrated circuit designers, but it would be favorable for the parallel algorithm researchers. There are four reasons below.

1. It isolates different circuits from each other, and one circuit does not need to know any details about other ones. This could be exactly explained by the Distributed Memory Access model.

2. It transfers the interfacial potentials and currents from one circuit to another, which could be considered as the message passing approach in parallel computing [8].

3. It only exists between adjacent circuits, so the communication just takes place between adjacent processors. This is an instance of the N2N communication model.

4. Its existence does not affect the stability of the resistor network. This observation is the physical base of the convergence theory of VTM.

Consequently, we may ask how to make use of the transmission line to boost the parallel computing of sparse linear systems. Obviously, there is no transmission line in this mathematical problem, so we have to add them artificially. VTM is then discovered. It inserts the Virtual Transmission Lines (VTL) into the linear system to achieve parallel computing.

## 3. Weighted Graph and Electric Graph.

In this section we define the weighted graph for the matrix, and define the electric graph for the linear system.

Assume there is an $n$-dimension linear system,

(3.1) $$\mathbf{Ax} = \mathbf{b}$$

where $\mathbf{x} = \begin{pmatrix} x_1 \\ x_2 \\ \vdots \\ x_n \end{pmatrix}$, $\mathbf{b} = \begin{pmatrix} b_1 \\ b_2 \\ \vdots \\ b_n \end{pmatrix}$, $\mathbf{A} = \begin{pmatrix} a_{11} & a_{12} & \cdots & a_{1n} \\ a_{21} & a_{22} & \cdots & a_{2n} \\ \vdots & \vdots & \ddots & \vdots \\ a_{n1} & a_{n2} & \cdots & a_{nn} \end{pmatrix}$, $a_{ij} = a_{ji}$, if $i \neq j$. $\mathbf{A}$ is symmetric.

As a symmetric matrix, $\mathbf{A}$ could be represented by an undirected graph $G$ [2, 4]. Each vertex $V_i$ of $G$ is one-to-one mapped to an unknown $x_i$ of the linear system. There is an edge $E_{ij}$ between $V_i$ and $V_j$ in $G$, iff $a_{ij} \neq 0$, $i \neq j$; otherwise, $V_i$ and $V_j$ are not connected.

A weighted graph $G_a$ is an undirected graph defined with the vertex weights and edge weights. $a_{ii}$ is defined as the weight of $V_i$, and $a_{ij}$, $i \neq j$, is defined as the weight of $E_{ij}$. A weighted graph is one-to-one mapped to a symmetric matrix. $G_a$ is defined to be SPD, if and only if the coefficient matrix $\mathbf{A}$ is SPD.

An electric graph $G_e$ is a weighted graph defined with the current sources. $b_i$ is defined as the inflow current source of $V_i$. We call $x_i$ the potential of $V_i$, and $\mathbf{x}$ is the potential vector of $G_e$. An electric graph is one-to-one mapped to a symmetric linear system. $G_e$ is defined to be SPD if and only if its corresponding $G_a$ is SPD.

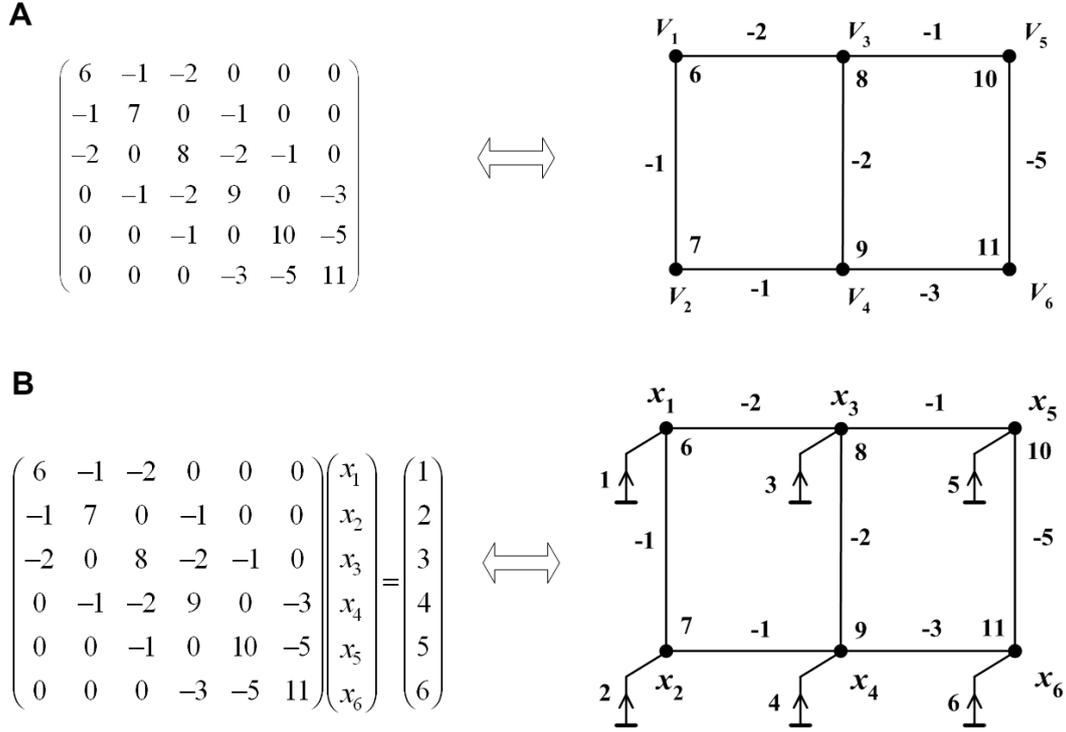

**Figure 3.** (A) The weighted graph of the matrix **A**. (B) The electric graph of the linear system **Ax** = **b**.

**Example 3.1**: The weighted graph of the coefficient matrix of (3.2) is shown in Fig. 3A, and the electric graph of this linear system is shown in Fig. 3B.

(3.2)
$$\begin{pmatrix} 6 & -1 & -2 & 0 & 0 & 0 \\ -1 & 7 & 0 & -1 & 0 & 0 \\ -2 & 0 & 8 & -2 & -1 & 0 \\ 0 & -1 & -2 & 9 & 0 & -3 \\ 0 & 0 & -1 & 0 & 10 & -5 \\ 0 & 0 & 0 & -3 & -5 & 11 \end{pmatrix} \begin{pmatrix} x_1 \\ x_2 \\ x_3 \\ x_4 \\ x_5 \\ x_6 \end{pmatrix} = \begin{pmatrix} 1 \\ 2 \\ 3 \\ 4 \\ 5 \\ 6 \end{pmatrix}$$

## 4. Electric Vertex Splitting

Before the parallel computing of the symmetric linear system **Ax** = **b**, we should partition it first. In this section, we introduce a new splitting technique to partition the electric graph of the symmetric linear system, which is called Electric Vertex Splitting (EVS). To partition the sparse linear system from circuit, EVS is also called wire tearing.

EVS is based on Kirchhoff's Current Law from electrical engineering [21]. The major difference of EVS over the traditional partitioning algorithm is that we bring in some new unknowns, called inflow currents, to the subgraphs. We may consider the electric graph to be a linear electric network, and we may recognize the vertex to be an electric node, and the edge to be a branch. An electric network has not only potentials but also currents. When one node is split into two twin vertices, the continuous current inside is also cut off and thus disclosed, so it is reasonable for us to consider these disclosed currents when doing the splitting. This concept is illustrated in Fig. 4.

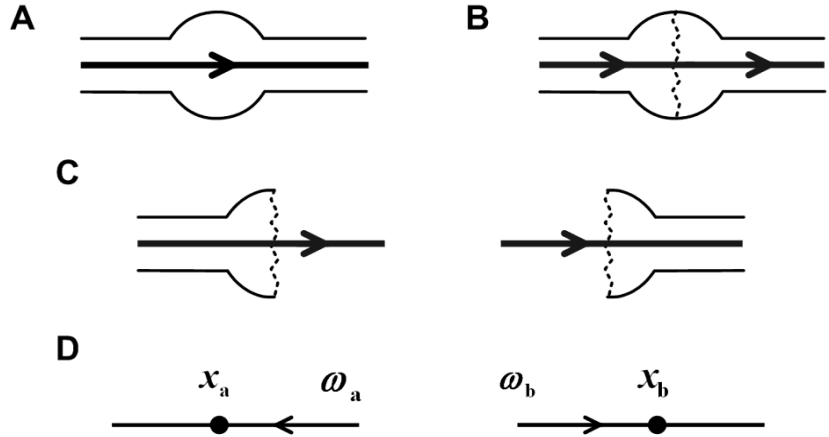

**Figure 4.** Illustration of the Electric Vertex Splitting. (**A**) The original node, with current flowing through it. (**B**) Splitting this node. (**C**) The node is split into a pair of twin nodes, and the currents are disclosed. (**D**) Simplified symbol of the inflow currents.

There are four steps to perform EVS upon the electric graph.

**Step 1.** Set the splitting boundary $G_B$. $V \in G_e$ is called boundary vertex iff $V \in G_B$; otherwise, $V$ is called inner vertex.

**Step 2.** Split each boundary vertex into a pair of vertices, which are called twin vertices. The original boundary vertex is called parent vertex.

**Step 3.** Split the weight and current source of each boundary vertex, and split the weight of each edge along the boundary, i.e. $E_{ij}$, if $E_{ij} \in G_e$ and $V_i, V_j \in G_B$.

**Step 4.** Add inflow currents to the twin vertices. These inflow currents represent the disclosed currents after splitting.

After these four steps, the original electric graph is split into $N$ subgraphs. If there is inflow current flowing into one vertex, then this vertex is called a port. As the result, twin vertices are also the ports of subgraghs.

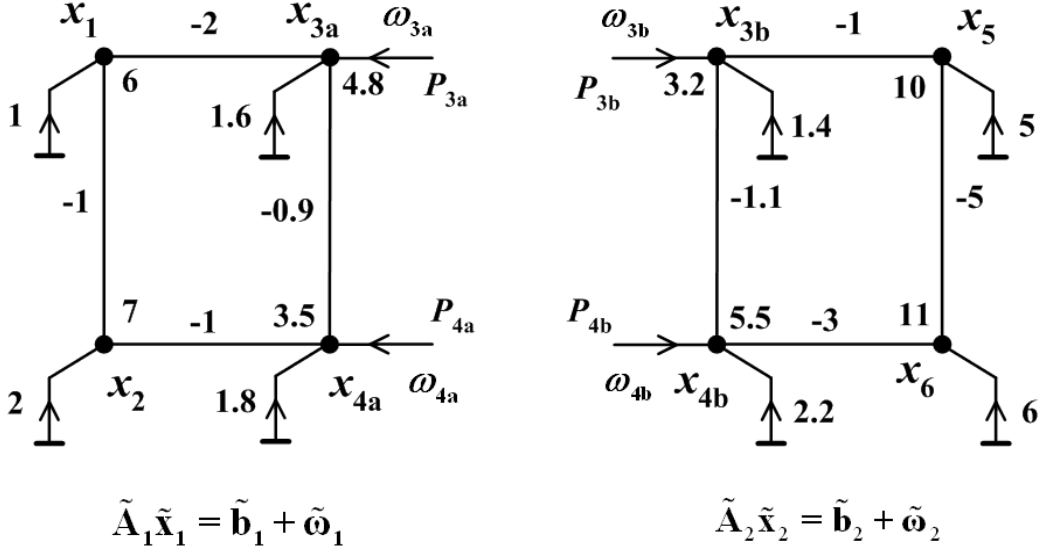

**Figure 5.** Electric Vertex Splitting of the electric graph of **Ax = b**

**Example 4.1**:

Continuing with Example 3.1, we split the electric graph $G_e$ of the linear system (3.2), previously shown in Fig. 3A. $V_3$ and $V_4$ are set to be the boundary $G_B$ and we split the weights and current sources of them. Please be noted that the weight of the edge $E_{34}$ is also split into two parts, $-0.9$ and $-1.1$. then we get 4 ports, $P_{3a}$, $P_{3b}$, $P_{4a}$ and $P_{4b}$, with currents $\omega_{3a}$, $\omega_{3b}$, $\omega_{4a}$ and $\omega_{4b}$ flowing into them, respectively. After that $G_e$ is split into two subgraghs. Finally we obtain two subsystems (4.1) and (4.2). Fig. 5 illustrates the process of EVS.

$$(4.1) \quad \begin{pmatrix} 6 & -1 & -2 & 0 \\ -1 & 7 & 0 & -1 \\ -2 & 0 & 4.8 & -0.9 \\ 0 & -1 & -0.9 & 3.5 \end{pmatrix} \begin{pmatrix} x_1 \\ x_2 \\ x_{3a} \\ x_{4a} \end{pmatrix} = \begin{pmatrix} 1 \\ 2 \\ 1.6 \\ 1.8 \end{pmatrix} + \begin{pmatrix} 0 \\ 0 \\ \omega_{3a} \\ \omega_{4a} \end{pmatrix}, \text{ or } \tilde{\mathbf{A}}_1 \tilde{\mathbf{x}}_1 = \tilde{\mathbf{b}}_1 + \tilde{\boldsymbol{\omega}}_1$$

$$(4.2) \quad \begin{pmatrix} 3.2 & -1.1 & -1 & 0 \\ -1.1 & 5.5 & 0 & -3 \\ -1 & 0 & 10 & -5 \\ 0 & -3 & -5 & 11 \end{pmatrix} \begin{pmatrix} x_{3b} \\ x_{4b} \\ x_5 \\ x_6 \end{pmatrix} = \begin{pmatrix} 1.4 \\ 2.2 \\ 5 \\ 6 \end{pmatrix} + \begin{pmatrix} \omega_{3b} \\ \omega_{4b} \\ 0 \\ 0 \end{pmatrix}, \text{ or } \tilde{\mathbf{A}}_2 \tilde{\mathbf{x}}_2 = \tilde{\mathbf{b}}_2 + \tilde{\boldsymbol{\omega}}_2$$

It should be noted that there are 12 unknowns in (4.1) and (4.2), while there are only 8 equations. Therefore, extra equations, also called boundary conditions, should be supplemented in order to construct an iterative relationship. Boundary conditions will be described in Theorem 4.1 and Section 5.

The split electric graph which consists of $N$ subgraghs is represented by $\tilde{G}_e$. Usually, there is more than one way to choose the splitting boundary, and even the splitting boundary is chosen, there are still plenty of ways to split the weights and current sources. Each of these ways is called a partition scheme of the electric graph.

EVS could also be used to split the weighted graph of a symmetric matrix **A**. Since no current sources in the weighted graph, it is unnecessary to add the currents into the twin vertex after splitting.

As the result, to split the weighted graph $G_a$ by EVS, there are three steps.

**Step 1.** Set the splitting boundary $G_B$, $G_B \subseteq G_a$.

**Step 2.** Split each boundary vertex $V \in G_B$ into a pair of twin vertices.

**Step 3.** Split the weight of each boundary vertex, and split the weight of each edge along the boundary, i.e. $E_{ij}$, if $E_{ij} \in G_a$ and $V_i, V_j \in G_B$.

**Example 4.2**:

Continuing with Example 4.1, we split the coefficient matrix **A** of linear system (3.2), whose weighted graph $G_a$ was previously shown in Fig. 3A.

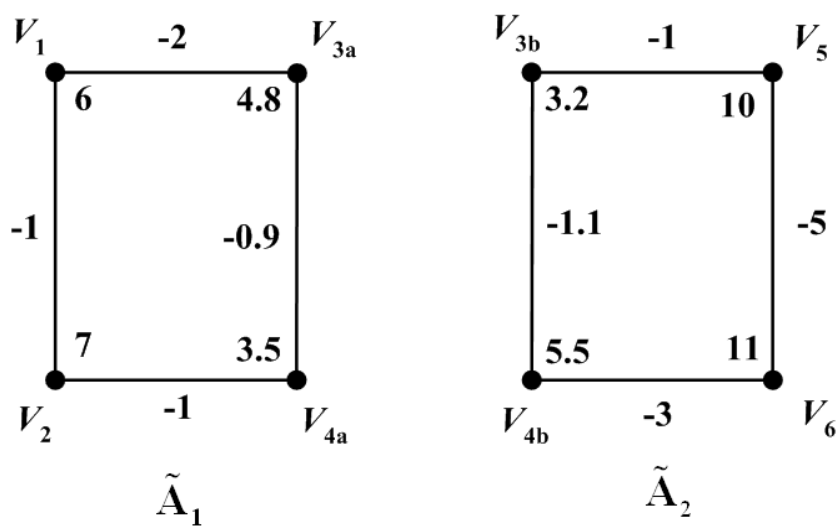

**Figure 6.** Electric vertex splitting of the weighted graph of the matrix **A**

We set $V_3$ and $V_4$ to be the boundary and split the weights of them. The split result is shown in Fig. 6. After that, $G_a$ is split into two subgraphs, which means that the original matrix **A** is split into two matrices, $\tilde{\mathbf{A}}_1$ and $\tilde{\mathbf{A}}_2$.

$$\tilde{\mathbf{A}}_1 = \begin{pmatrix} 6 & -1 & -2 & 0 \\ -1 & 7 & 0 & -1 \\ -2 & 0 & 4.8 & -0.9 \\ 0 & -1 & -0.9 & 3.5 \end{pmatrix}, \quad \tilde{\mathbf{A}}_2 = \begin{pmatrix} 3.2 & -1.1 & -1 & 0 \\ -1.1 & 5.5 & 0 & -3 \\ -1 & 0 & 10 & -5 \\ 0 & -3 & -5 & 11 \end{pmatrix}$$

After illustrating an example of EVS, we present its mathematical description. Assume the original graph $G_e$ is partitioned into $N$ separated subgraghs, $M_j, j = 1, 2, \cdots, N$, following some partition scheme. Thereafter, we use $|M_j|$ to represent the number of vertices in $M_j$. $M_j$ and $M_i$ are called adjacent subgraghs, if each of them has at least one twin vertex born from the same boundary vertex.

Each subgragh could be mapped back into a symmetric linear subsystem with inflow currents. To express this subsystem, we define $\Gamma_{j,port}$ to be an ordered set of the ports in $M_j$, and $\Gamma_{j,inner}$ an ordered set of the inner vertices in $M_j$. We define $\mathbf{u}_j$ to be the potential vector of $\Gamma_{j,port}$, and $\mathbf{y}_j$ to be the potential vector of $\Gamma_{j,inner}$. Then, the local linear system for each subgragh could be expressed by the following equation:

$$(4.3) \quad \begin{bmatrix} \mathbf{C}_j & \mathbf{E}_j \\ \mathbf{F}_j & \mathbf{D}_j \end{bmatrix} \begin{bmatrix} \mathbf{u}_j \\ \mathbf{y}_j \end{bmatrix} = \begin{bmatrix} \mathbf{f}_j \\ \mathbf{g}_j \end{bmatrix} + \begin{bmatrix} \boldsymbol{\omega}_j \\ 0 \end{bmatrix}$$

where $j = 1, 2, \cdots, N$. $\boldsymbol{\omega}_j$ is the inflow current vector of the ports of $M_j$. The inflow current of an inner vertex is zero. $\mathbf{u}_j$ and $\boldsymbol{\omega}_j$ are also called the local boundary variables of $M_j$, respectively.

The above equations (4.3) could be simply rewritten as,

(4.4) $$\tilde{\mathbf{A}}_j \tilde{\mathbf{x}}_j = \tilde{\mathbf{b}}_j + \tilde{\boldsymbol{\omega}}_j$$

where $j = 1, 2, \cdots, N$.

The above-mentioned splitting technique is called level-one splitting technique, and the split vertices could be split again and again, which are called multilevel splitting technique, as illustrated in Fig. 7. To partition a physical problem in 2 or 3 dimensions, the level-two and level-three splitting techniques are inevitable.

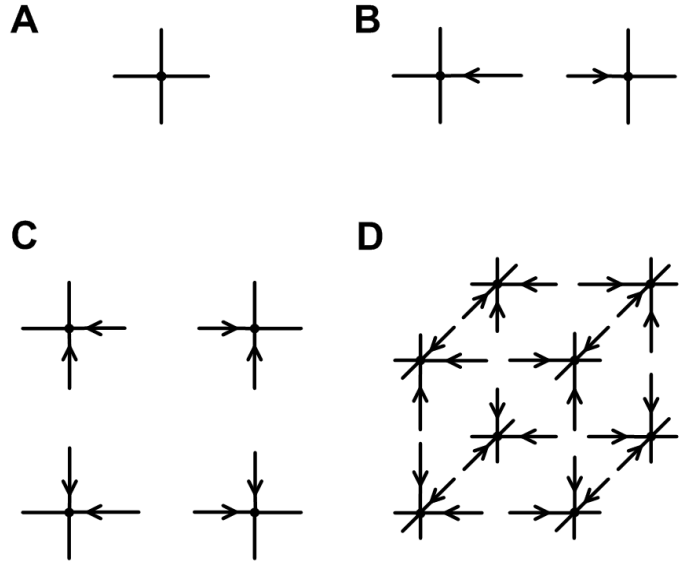

**Figure 7.** Illustration of multilevel Electric Vertex Splitting.

(**A**) The original vertex.

(**B**) Level-one splitting, where one vertex is split into a pair of twin vertices, and there is one inflow current into each twin vertex.

(**C**) Level-two splitting, where one vertex is split into four child vertices, and there are two inflow currents into each child vertex.

(**D**) Level-three splitting, where one vertex is split into eight child vertices, and there are three inflow currents into each child vertex.

**Theorem 4.1** (Reversibility): Suppose the electric graph $G_e$ is partitioned into $\tilde{G}_e$ by Electric Vertex Splitting (EVS). If the potentials of each pair of twin vertices are set to be same, and the inflow currents of them are set to be opposite, then $\tilde{G}_e$ is equivalent to $G_e$, i.e. the potential of each pair of twin vertices in $\tilde{G}_e$ is equal to the potential of their parent vertex in $G_e$, and the potential of each inner vertex in $\tilde{G}_e$ is equal to that of its corresponding inner vertex in $G_e$.

This theorem tells us that EVS is reversible, and this is easy to understand according to its physical background. If we reverse the process of EVS, which means that we make the inflow currents to be a continuous current, merge the twin vertices into one vertex and envelop the continuous current inside it, then we get the original electric graph. A proof for this theorem is given in Appendix 2.

**Example 4.3**:

Continuing with Example 4.1, we set:

(4.5)
$$\begin{cases} x_{3a} = x_{3b} = x_3 \\ x_{4a} = x_{4b} = x_4 \\ \omega_{3a} + \omega_{3b} = 0 \\ \omega_{4a} + \omega_{4b} = 0 \end{cases}$$

Combining (4.1), (4.2) and (4.5), we get (3.2) after eliminating $x_{3a}$, $x_{3b}$, $\omega_{3a}$, $\omega_{3b}$, $x_{4a}$, $x_{4b}$, $\omega_{4a}$ and $\omega_{4b}$.

**Theorem 4.2** (Conformal Splitting Existence): Suppose the weighted graph $G_a$ is SPD, then for arbitrarily-chosen boundary, there is more than one scheme to partition $G_a$ into $N$ subgraphs by EVS, $\tilde{G}_j, j = 1, 2, \cdots, N$, and all $\tilde{G}_j$ are SPD.

This theorem assures that an SPD graph must be able to be partitioned into arbitrary number of SPD subgraphs by EVS. A proof is given in Appendix 1. Here we reuse the word "conformal" to represent a kind of EVS partition schemes, which hold the SPD property of the electric graph.

For the electric graph, we have the same conclusion, as below:

Suppose the electric graph $G_e$ is SPD, then for arbitrarily-chosen boundary, there is more than one scheme to partition $G_e$ into $N$ subgraghs by EVS, $\tilde{G}_j, j = 1, 2, \cdots, N$, and all $\tilde{G}_j$ are SPD.

**Corollary 4.1:** Suppose the electric graph $G_e$ is SPD, then for arbitrarily-chosen boundary, there is more than one scheme to partition $G_e$ into $N$ subgraghs by EVS, $\tilde{G}_j, j = 1, 2, \cdots, N$, which are all symmetric-nonnegative-definite (SNND).

Corollary 4.1 is a weak version of Theorem 4.2. Then, we broaden the definition of conformal partition:

If some partition scheme makes Corollary 4.1 work, then this scheme is conformal, since it holds the SPD or SNND property for the subgraghs after partitioning.

This paper does not figure out how to set a practical partition scheme for EVS to split the electric graph conformally. This is a simple work for any strongly-diagonal or weakly-diagonal sparse system. For the scientific problem, we recommend to do the partitioning on the physical level before generating sparse linear systems.

## 5. VTM

Assume that the electric graph $G_e$ has been partitioned into $N$ subgraghs, then we add one VTL between each pair of twin vertices, which means that we use the transmission equations as the boundary conditions. A simple example is given as below.

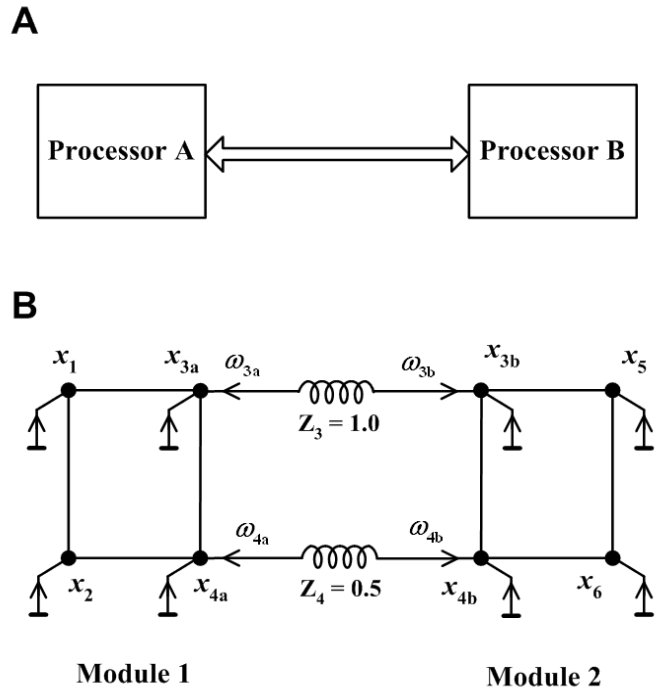

**Figure 8.** The split electric graph with VTLs.

**Example 5.1**: Continuing with Example 4.1, we add one VTL $T_3$ between $x_{3a}$ and $x_{3b}$, whose characteristic impedance $Z_3$ is set to be 1, then we add another line $T_4$ between $x_{4a}$ and $x_{4b}$, whose $Z_4$ is set to be 0.5.

According to (2.1), the mathematical equation of $T_3$ is:

(5.1)
$$\begin{cases} x_{3a}^{k} + 1.0 \cdot \omega_{3a}^{k} = x_{3b}^{k-1} - 1.0 \cdot \omega_{3b}^{k-1} \\ x_{3b}^{k} + 1.0 \cdot \omega_{3b}^{k} = x_{3a}^{k-1} - 1.0 \cdot \omega_{3a}^{k-1} \end{cases}$$

where $k$ is the iteration index in VTM.

Similarly, the mathematical equation of $T_4$ is:

(5.2)
$$\begin{cases} x_{4a}^{k} + 0.5 \cdot \omega_{4a}^{k} = x_{4b}^{k-1} - 0.5 \cdot \omega_{4b}^{k-1} \\ x_{4b}^{k} + 0.5 \cdot \omega_{4b}^{k} = x_{4a}^{k-1} - 0.5 \cdot \omega_{4a}^{k-1} \end{cases}$$

Based on (4.1) and part of (5.1) and (5.2), the linear system of Subgragh 1 could be expressed as below:

(5.3)
$$\begin{cases} \begin{pmatrix} 6 & -1 & -2 & 0 \\ -1 & 7 & 0 & -1 \\ -2 & 0 & 4.8 & -0.9 \\ 0 & -1 & -0.9 & 3.5 \end{pmatrix} \begin{pmatrix} x_1 \\ x_2 \\ x_{3a} \\ x_{4a} \end{pmatrix} = \begin{pmatrix} 1 \\ 2 \\ 1.6 \\ 1.8 \end{pmatrix} + \begin{pmatrix} 0 \\ 0 \\ \omega_{3a} \\ \omega_{4a} \end{pmatrix} \\ x_{3a}^{k} + 1.0 \cdot \omega_{3a}^{k} = x_{3b}^{k-1} - 1.0 \cdot \omega_{3b}^{k-1} \\ x_{4a}^{k} + 0.5 \cdot \omega_{4a}^{k} = x_{4b}^{k-1} - 0.5 \cdot \omega_{4b}^{k-1} \end{cases}$$

Eliminate $\omega_{3a}^{k}$ and $\omega_{4a}^{k}$ from (5.3) and we get (5.4):

(5.4)
$$\begin{pmatrix} 6 & -1 & -2 & 0 \\ -1 & 7 & 0 & -1 \\ -2 & 0 & 5.8 & -0.9 \\ 0 & -1 & -0.9 & 5.5 \end{pmatrix} \begin{pmatrix} x_1 \\ x_2 \\ x_{3a} \\ x_{4a} \end{pmatrix} = \begin{pmatrix} 1 \\ 2 \\ 1.6 + x_{3b}^{k-1} - \omega_{3b}^{k-1} \\ 1.8 + 2 \cdot x_{4b}^{k-1} - \omega_{4b}^{k-1} \end{pmatrix}$$

$$\begin{cases} \omega_{3a}^{k} = -x_{3a}^{k} + x_{3b}^{k-1} - \omega_{3b}^{k-1} \\ \omega_{4a}^{k} = -2x_{4a}^{k} + 2x_{4b}^{k-1} - \omega_{4b}^{k-1} \end{cases}$$

Similarly, we get (5.5) for Subgragh 2:

$$(5.5) \quad \begin{pmatrix} 4.2 & -1.1 & -1 & 0 \\ -1.1 & 7.5 & 0 & -3 \\ -1 & 0 & 10 & -5 \\ 0 & -3 & -5 & 11 \end{pmatrix} \begin{pmatrix} x_{3b} \\ x_{4b} \\ x_5 \\ x_6 \end{pmatrix} = \begin{pmatrix} 1.4 + x_{3a}^{k-1} - \omega_{3a}^{k-1} \\ 2.2 + 2 \cdot x_{4a}^{k-1} - \omega_{4a}^{k-1} \\ 5 \\ 6 \end{pmatrix}$$

$$\begin{cases} \omega_{3b}^k = -x_{3b}^k + x_{3a}^{k-1} - \omega_{3a}^{k-1} \\ \omega_{4b}^k = -2x_{4b}^k + 2x_{4a}^{k-1} - \omega_{4a}^{k-1} \end{cases}$$

After that, we set the initial value of the boundary variables as below:

$$\begin{cases} x_{3a}^0 = x_{3b}^0 = x_{4a}^0 = x_{4b}^0 = 0 \\ \omega_{3a}^0 = \omega_{3b}^0 = \omega_{4a}^0 = \omega_{4b}^0 = 0 \end{cases}$$

At last, we compute this example distributedly on two processors. Subgragh 1 is located on Processor A, and Subgragh 2 is located on Processor B, as illustrated in Fig. 8. The boundary variables are communicated between these two processors by message passing. The computing result is shown in Fig. 9.

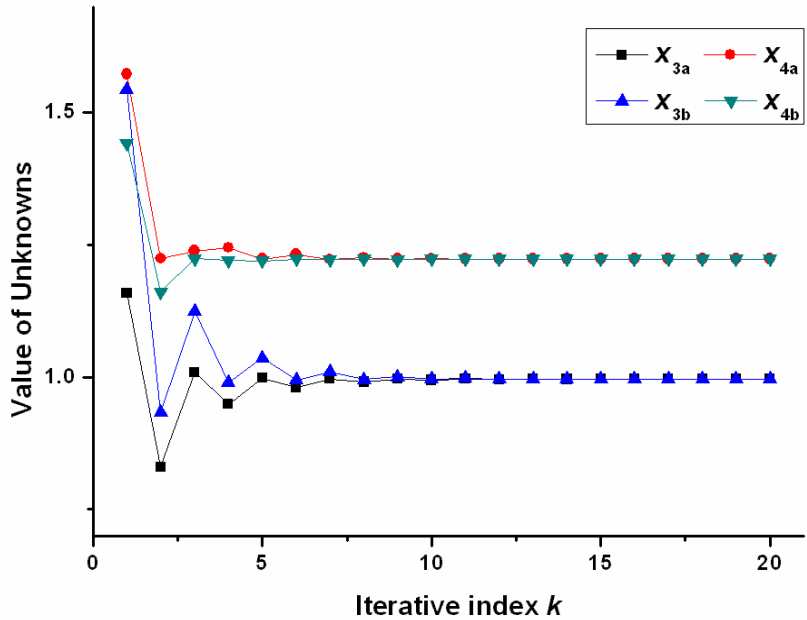

**Figure 9.** Distributed computing result of VTM on double processors

After illustrating this simple example, we present the mathematical description of VTM. For the subgragh $M_j$, we have defined $\Gamma_{j,port}$ as an ordered set of its ports. Further, we define $\Gamma_{j,twin}$ to be another ordered set of ports whose twin vertices belong to $\Gamma_{j,port}$. The ports in $\Gamma_{j,port}$ and their corresponding twin vertices in $\Gamma_{j,twin}$ have the same order. It is easy to know that the ports of $\Gamma_{j,twin}$ belong to the adjacent subgraghs of $M_j$. Define $\mathbf{u}_{j,twin}$ as the potential vector of $\Gamma_{j,twin}$, and $\boldsymbol{\omega}_{j,twin}$ as the current vector of $\Gamma_{j,twin}$. Then, for $M_j$, $j = 1, 2, \cdots, N$, the transmission equations (2.1) could be expressed in a matrix-vector form as below:

$$\mathbf{u}_j^k + \mathbf{Z}_j \boldsymbol{\omega}_j^k = \mathbf{u}_{j,twin}^{k-1} - \mathbf{Z}_j \boldsymbol{\omega}_{j,twin}^{k-1} \tag{5.6}$$

Here $\mathbf{Z}_j$ is a positive diagonal matrix, called the local characteristic impedance matrix of $M_j$. The diagonal elements of $\mathbf{Z}_j$ are the characteristic impedances of the VTLs connected to $M_j$. $\mathbf{Z}_j$ is the local preconditioner for subgraph $\mathbf{M}_j$.

(5.6) is an distributedly-iterative relation, and $\mathbf{u}_{j,twin}^{k-1}$ and $\boldsymbol{\omega}_{j,twin}^{k-1}$ are the previous computing results passed from the adjacent subgraghs, which are called the remote boundary variables of $M_j$. Merge (4.1) and (5.6), we get:

$$\begin{bmatrix} \mathbf{C}_j & \mathbf{E}_j & -\mathbf{I} \\ \mathbf{F}_j & \mathbf{D}_j & 0 \\ \mathbf{I} & 0 & \mathbf{Z}_j \end{bmatrix} \begin{bmatrix} \mathbf{u}_j^k \\ \mathbf{y}_j^k \\ \boldsymbol{\omega}_j^k \end{bmatrix} = \begin{bmatrix} \mathbf{f}_j \\ \mathbf{g}_j \\ \mathbf{u}_{j,twin}^{k-1} - \mathbf{Z}_j \boldsymbol{\omega}_{j,twin}^{k-1} \end{bmatrix} \tag{5.7}$$

where $\mathbf{I}$ is the identity matrix. Eliminating $\boldsymbol{\omega}_j^k$, we get the following SPD system:

$$\begin{bmatrix} \mathbf{C}_j + \mathbf{Z}_j^{-1} & \mathbf{E}_j \\ \mathbf{F}_j & \mathbf{D}_j \end{bmatrix} \begin{bmatrix} \mathbf{u}_j^k \\ \mathbf{y}_j^k \end{bmatrix} = \begin{bmatrix} \mathbf{f}_j + \mathbf{Z}_j^{-1}(\mathbf{u}_{j,twin}^{k-1} - \mathbf{Z}_j \boldsymbol{\omega}_{j,twin}^{k-1}) \\ \mathbf{g}_j \end{bmatrix} \tag{5.8}$$

$$\boldsymbol{\omega}_j^k = -\mathbf{Z}^{-1} \cdot \mathbf{u}_j^k + \mathbf{Z}^{-1} \cdot \mathbf{u}_{j,twin}^{k-1} - \boldsymbol{\omega}_{j,twin}^{k-1}$$

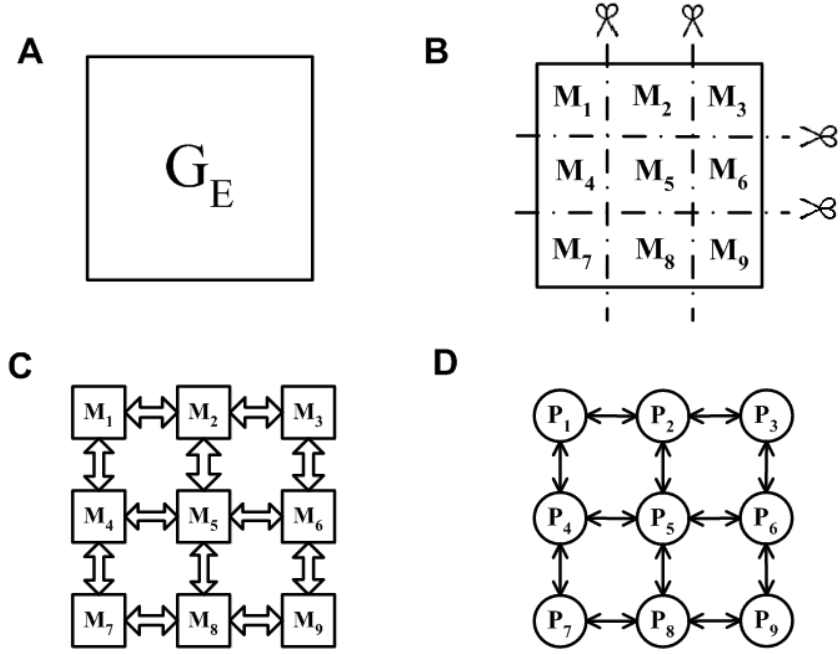

**Figure 10.** Illustration of the computing process of VTM.

(**A**) The original electric graph of the sparse linear system.

(**B**) Partition the original graph into $N$ subgraghs by EVS.

(**C**) Add VTLs between adjacent subgraghs.

(**D**) Map each subgragh onto one processor.

(5.8) is called the local subsystem of $M_j$, which could be solved by Sparse or Dense Cholesky, CG, MG, etc.

Table 1 gives the full description of VTM, and Fig. 10 illustrates the computing process of this algorithm. It should be noted that there is no broadcasting, but only N2N communication.

**Table 1.** Algorithm description of VTM

---

Assume the original electric graph has been partitioned into $N$ subgraghs. Each subgragh is located on one processor, and there are communication networks between adjacent subgraghs.

For Subgragh $j$, $j = 1, \cdots, N$, do in parallel:

1. Communicate with adjacent subgraghs, to make an agreement of the characteristic impedances for each VTL, so that $\mathbf{Z}_j$ is set.

2. Guess the initial local boundary variable $\mathbf{u}_j^0$ and $\boldsymbol{\omega}_j^0$ of each port, which could be arbitrary values.

3. Wait until receiving the new remote boundary variables, $\mathbf{u}_{j,twin}^{k-1}$ and $\boldsymbol{\omega}_{j,twin}^{k-1}$, do:

4. Solve the local subsystem with the updated remote boundary variable $\mathbf{u}_{j,twin}^{k-1}$ and $\boldsymbol{\omega}_{j,twin}^{k-1}$, and then we get the new local boundary variable $\mathbf{u}_j^k$ and $\boldsymbol{\omega}_j^k$.

5. Send the new local boundary variable $\mathbf{u}_j^k$ and $\boldsymbol{\omega}_j^k$ to the adjacent subgraghs.

6. If convergent,
    Break;

7. EndWait

## 6. Convergence theory of VTM.

According to the description of VTM, it is not straightforward to judge whether this algorithm is convergent or not. In this section we present the convergence theorem.

**Theorem 6.1** (Convergence): Assume the electric graph of an SPD linear system, $\mathbf{Ax} = \mathbf{b}$, is partitioned into $N$ symmetric-non-negative-definite (SNND) subgraghs, then for positive characteristic impedances of VTLs, VTM converges to the solution of the original system.

This conclusion is valid for both the level-one and the multilevel EVS, and we give a proof for this theorem in Appendix 3.

Theorem 6.1 could also be simplified as: Assume the electric graph of a SPD linear system is partitioned into $N$ subgraghs following a conformal partition scheme, then VTM converges.

## 7. Preconditioning

As we observed, the choice of the characteristic impedances of VTLs, would make a huge impact to the convergence speed of VTM. Consequently, the characteristic impedances, i.e. the characteristic impedance matrix $\mathbf{Z}_j$, could be considered as the preconditioner for VTM. Further, we define the preconditioning of VTM as the process to find proper characteristic impedance matrix for VTLs.

## 7.1. Impedance Matching

Here we propose a simple way, called impedance matching, to choose the characteristic impedances, i.e. to precondition VTM.

Before describing this technique, it is necessary to define the port's input impedance, which could be found in any textbook of circuit theory or microwave network.

The theory of VTM could be considered as a mix of numerical analysis and microwave network. This paper borrows quite a few notations and definitions from electrical engineering, such as transmission line, potential, source current, inflow current, characteristic impedance, etc.

**Definition 7.1** (Input Impedance of Port):

For the subgragh described by (4.3), we first set all the inflow current sources to be zero, and then set the inflow currents of all the ports except $P_j$ to be zero, and set the inflow current $\omega_j$ of $P_j$ to be 1, than we solve this system and get the potential $u_j$ of $P_j$, then, $r_{in,j} = u_j / \omega_j = u_j / 1 = u_j$, here $r_{in,j}$ is the input impedance of port $P_j$.

The impedance matching technique is that, the characteristic impedance of VTL should be neither too large nor too small, and usually it is set near the input impedances of either port of VTL. We use the following example to illustrate the effect of impedance matching.

**Example 7.1**: We continue to use Example 5.1. The input impedance of $P_{3a}$ should be the answer of $x_{3a}$ in (7.3.1), and we get $r_{in,3a} = 0.2598$.

$$(7.3.1) \quad \begin{pmatrix} 6 & -1 & -2 & 0 \\ -1 & 7 & 0 & -1 \\ -2 & 0 & 4.8 & -0.9 \\ 0 & -1 & -0.9 & 3.5 \end{pmatrix} \begin{pmatrix} x_1 \\ x_2 \\ x_{3a} \\ x_{4a} \end{pmatrix} = \begin{pmatrix} 0 \\ 0 \\ 1 \\ 0 \end{pmatrix}$$

Similarly, we get $r_{in,4a} = 0.3190$, $r_{in,3b} = 0.3699$ and $r_{in,4b} = 0.2557$.

After that, we choose different combination of $Z_3$ and $Z_4$, and redo the computation in Example 5.1. The root mean squared (RMS) errors after 20 iterations are shown in Fig. 11, from which we know that the computational error of VTM is lowest when $Z_2$ is set near $r_{in,3a}$ or $r_{in,3b}$, and $Z_3$ is set near $r_{in,4a}$ or $r_{in,4b}$. This simple example shows that impedance matching is impactful to make VTM accurate and fast.

Then we test VTM on 128 processors and Fig. 12 illustrates the convergence curves of VTM with and without impedance matching, which is also impressive.

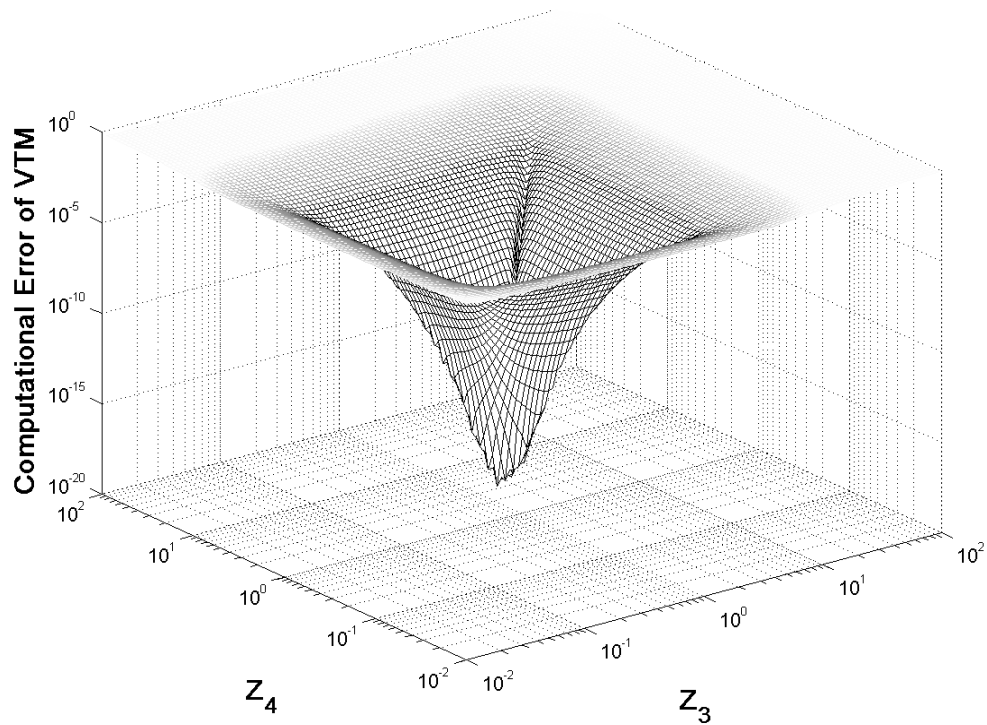

**Figure 11.** Computational error of VTM after 20 iterations

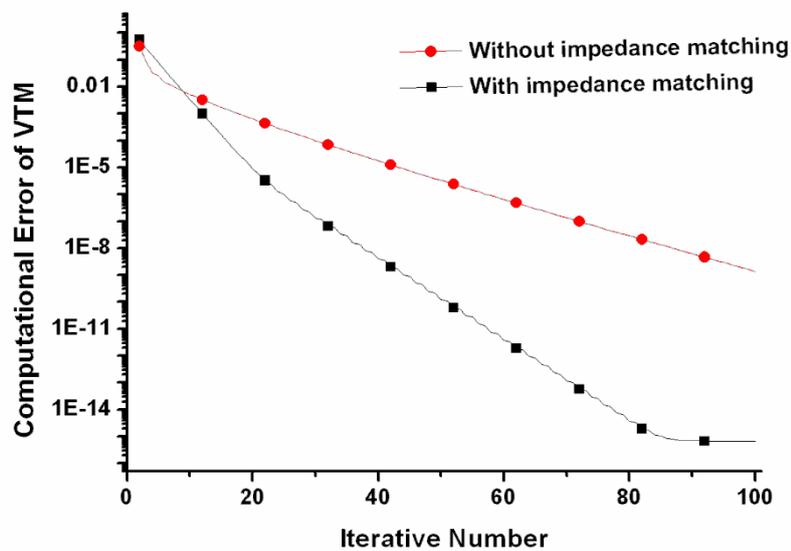

**Figure 12.** Effect of the impedance matching technique on 128 processors

At last, it should be noted that the computational error of VTM is a continuous function of the characteristic impedances of VTL, and it is not sensitive to the small change of the characteristic impedances. This character makes VTM to be a practical and robust numerical algorithm.

### 7.2. Coupling

Generally, the local preconditoner $\mathbf{Z}_j$ in (5.8) could not only be a diagonal matrix, but also a banded matrix or even a full matrix. In this case, these exists coupling among the adjacent VTLs. According to the knowledge of microwave network, if $\mathbf{Z}_j$ is a symmetrical matrix, then the VTLs connected to $M_j$ are symmetric coupled; if $\mathbf{Z}_j$ is an unsymmetrical matrix, these VTLs are unsymmetrical coupled. If $\mathbf{Z}_j$ is diagonal, the VTLs are uncoupled.

The microwave network with symmetric coupled transmission lines inclines to be more stable than that with uncoupled transmission lines. This means that the convergence of VTM with coupled VTLs might be faster than that with uncoupled VTLs.

If there exist coupled VTLs, the convergence theory of VTM is updated as below:

**Theorem 7.1**: Assume the electric graph of an SPD linear system, $\mathbf{Ax} = \mathbf{b}$, is partitioned into $N$ symmetric-non-negative-definite (SNND) subgraghs. If all the local preconditioner $\mathbf{Z}_j$ is SPD, VTM converges to the solution of the original system.

The proof for this theorem is similar to Theorem 6.1.

## 8. Performance Modeling

In this section we set a simple model for VTM [1, 22]. First we make several assumptions.

(1). One floating point operation at top speed (i.e. the speed of matrix multiplication) costs one time unit.

(2). We have $p$ processor arranged in a 2D mesh.

(3). The communication delay between neighboring processors are same, and sending a message of $l$ words from one processor to its adjacent processor costs ($\alpha + \beta * l$) time units.

Second, we prepare the linear system for test. The electric graph of this linear system is a 2D grid, whose dimension is $n$. By EVS, we partition this graph regularly into $p$ subgraghs. Each subgragh is a smaller grid whose dimension $b = \frac{n}{p} + 2\sqrt{\frac{n}{p}}$, and $b \approx \dfrac{n}{p}$.

Third, we locate each subgragh on one processor and use the sparse Cholesky factorization to solve the local subsystem. Numerical experiments shows that, to solve this kind of sparse linear systems, the computational complexity of sparse Cholesky factorization is $O(b^{1.5})$, and the computational complexity of the forward and backward substitution is $O(b)$.

Fourth, we do the precondition for VTM. The characteristic impedances of VTLs are optimized by impedance matching technique. The computational complexity of impedance matching is $O(b)$.

Fifth, we do the distributed iterative computation using VTM. Assume it needs $K$ iterations to achieve the computational error of $\varepsilon$. We need to do the Cholesky factorization for one time, and do the forward and backward substitution for the rest $K$-1 times, as explained in Section 7. Then, the total parallel computing time is:

(8.1)
$$T_p = b^{1.5} + K\left(2b + \alpha + \beta\sqrt{b}\right)$$
$$\approx \left(\frac{n}{p}\right)^{1.5} + 2K\left(\frac{n}{p}\right) + K\alpha + K\beta\left(\frac{n}{p}\right)^{0.5}$$

Compared to the computing time on a single processor:

(8.2)
$$T_s = n^{1.5} + 2n$$

The speedup ratio is:

$$S = \frac{T_s}{T_p} = \frac{n^{1.5} + 2n}{\left(\dfrac{n}{p}\right)^{1.5} + 2K\left(\dfrac{n}{p}\right) + K\alpha + K\beta\left(\dfrac{n}{p}\right)^{0.5}}$$
$$\approx \frac{1}{1 + 2K\sqrt{\dfrac{p}{n}} + K\beta\left(\dfrac{p}{n}\right)} \cdot p^{1.5}$$

Here the key is to know the total iterative number $K$, which could be approximately considered as a function of $n$ and $p$, i.e. $K(n, p)$. It is difficult to make a theoretical analysis of $K(n, p)$; however, numerical experiments in Section 9 show that the convergent speed of VTM is acceptable and $K$ is a moderate number to achieve high computational accuracy.

## 9. Numerical Experiments.

We test VTM by the VTM toolbox, which is a distributed computing emulation platform developed by us under MATLAB and SIMULINK. Here $n$ is the dimension of the sparse linear system, and $p$ is the number of cores.

We first test a sparse linear system whose dimension $n$ is 4225. We partition it into $p$ subgraghs and solve it on $p$ processors. Fig. 13 illustrates the RMS errors' curve of VTM when $p$ is 4, 8, 16, 32, 64 and 128. According to this figure, we know that the computational error of VTM is decreasing, and it is limited by the machine precision of the computer, which is double-precision in this case.

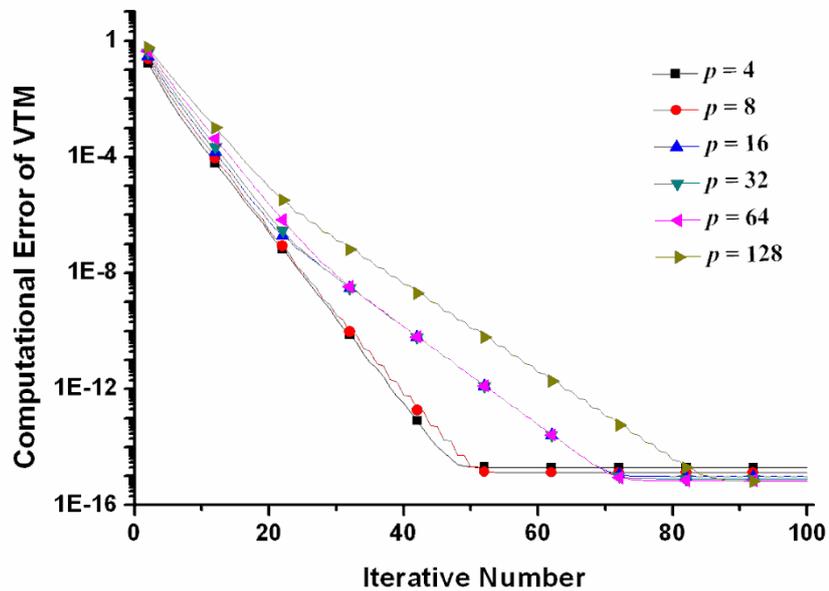

**Figure 13.** Computational errors of VTM when $p$ changes

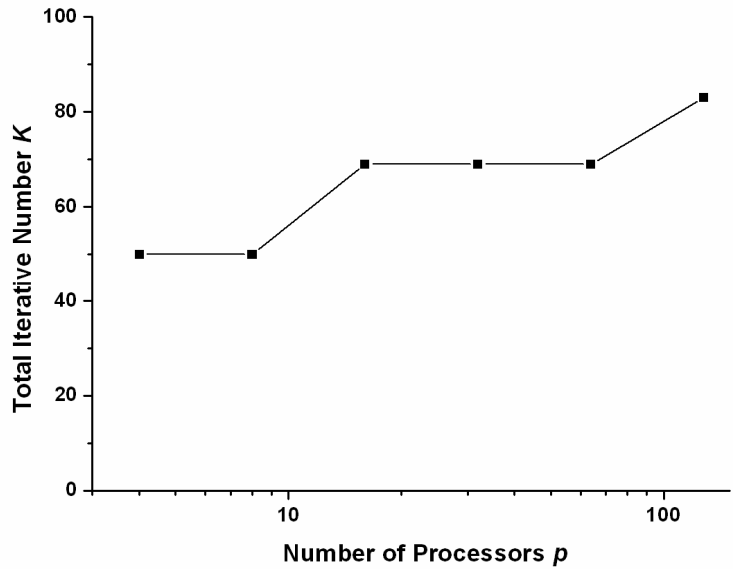

**Figure 14.** Illustration of *K(p)* when $\varepsilon = 2.0E-15$

If we set $\varepsilon = 2.0E-15$, then we get *K*(*p*) in Fig. 14, which is based on Fig. 13. This figure indicates that *K* increases slowly with *p*.

Then we solve a number of sparse linear systems on 64 processors. The dimensions of these testbenches are 289, 1089, 2401, 9409 and 14641, respectively. Fig. 15 illustrates the RMS errors' curve of VTM depending on *n*.

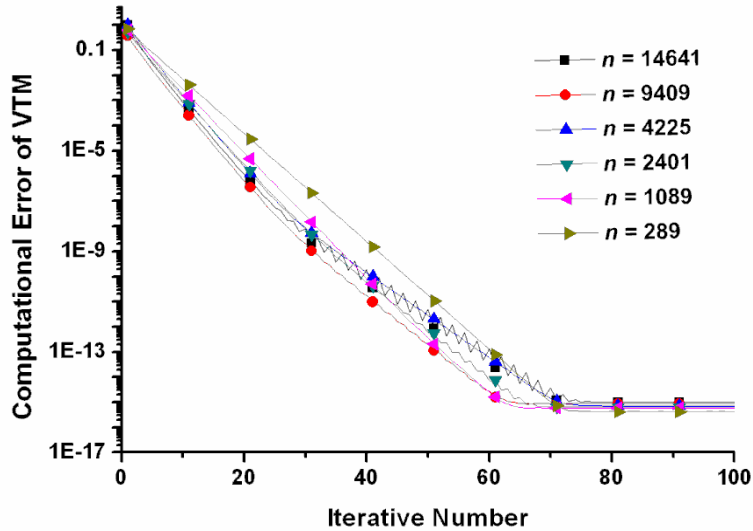

**Figure 15.** Computational errors of VTM when *n* changes, *p* = 64

If we set $\varepsilon = 2.0E-15$, then we get *K*(*n*) in Fig. 16. This figure indicates that *K* is somewhat immune to the change of *n*.

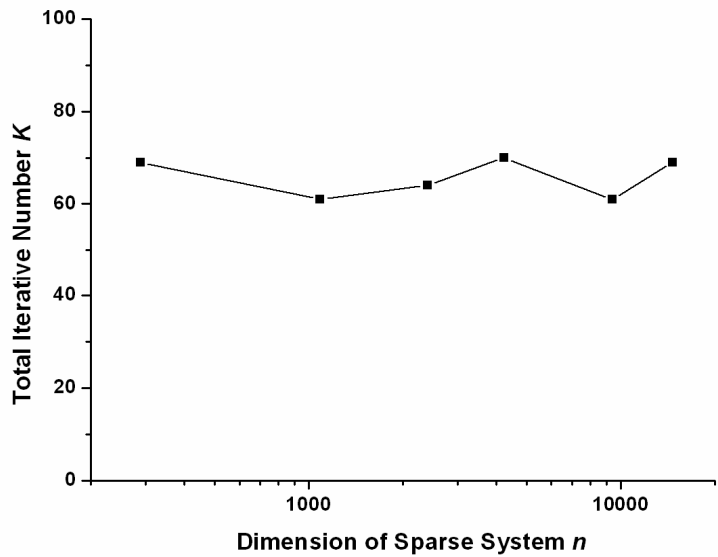

**Figure 16.** Illustration of *K(n)* when $\varepsilon = 2.0E-15$

These experiments show that VTM is an efficient and accurate algorithm. The total iteration number *K* is not sensitive to the change of *n*, which is the dimension of the sparse system, and *K* increases slowly with the number of processors *p*. As the result, if the dimension of subsystem on each processor were large enough, the efficiency of VTM might approach *p*, as predicted in Section 8.

Theoretically, the dimension of the sparse linear system being solved by VTM could be arbitrarily-large, and the processors being employed could be arbitrary number. Limited by our hardware, we are not able to test extremely large problem on supercomputers.

## 10. Conclusion

In this paper, we propose a new parallel algorithm, VTM, to solve the sparse SPD linear systems. VTM could be considered as a new block relaxation method, similar to block Jacobi, or a new algebraic domain decomposition algorithm, similar to additive Schwarz method.

The partitioning technique for VTM, i.e. Electric Vertex Splitting, is different from the traditional decomposition algorithms for sparse linear system.

The preconditioning of VTM is flexible. The characteristic impedance matrix has a strong impact to the convergence speed. If there is coupling between adjacent VTLs, the precondioner might be more efficient.

VTM could not only be used to solve the SPD systems, but also the non-SPD, unsymmetrical linear systems and nonlinear systems. For the unsymmetrical linear system, coupling technique would be helpful to make the algorithm easier to converge.

Acknowledgments.

We thank Prof. Hao Zhang, Yi Su, Dr. Chun Xia, Dr. Wei Xue, Pei Yang, Bin Niu. This work was partially sponsored by the Major State Basic Research Development Program of China (973 Program) under contract G1999032903, the National Natural Science Foundation of China Key Program, 90207001, and the National Science Fund for Distinguished Young Scholars of China, 60025101.

## Appendix 1. Proof of the Conformal Splitting Existence Theorem.

Before proving the Conformal Splitting Existence Theorem (Theorem 4.2), first we prove its simple version as below:

**Lemma A1.1**: Suppose the weighted graph $G_a$ is SPD, then for arbitrarily-chosen boundary $G_B$, there is more than one partition scheme to partition $G_a$ into two SPD subgraphs.

In order to prove Lemma A1.1, we present three other lemmas.

**Lemma A1.2**: The symmetric matrix $\mathbf{A}$ is one to one mapped to the quadratic form $P_\mathbf{A}(\mathbf{x}) = \mathbf{x}^T \mathbf{A} \mathbf{x}$.

**Lemma A1.3**: $\mathbf{A}$ is SPD, if and only if $P_\mathbf{A}(\mathbf{x})$ is positive-definite, i.e. $P_\mathbf{A}(\mathbf{x}) = \mathbf{x}^T \mathbf{A} \mathbf{x} > 0$, $\forall \mathbf{x} \in \mathbb{R}^n$, $\mathbf{x} \neq \mathbf{0}$.

According to Lemma A1.2 and A1.3, we know that Lemma A1.4 is right.

**Lemma A1.4**: To partition a weighted graph using the electric vertex splitting is equivalent to divide its quadratic form using the variable splitting technique, and vice versa.

The following example illustrates the variable splitting technique for the quadratic form.

**Example A1.1**: This example is based on Example 4.1 in Section 4. The quadratic form of $\mathbf{A}$ is:

(A1.1) $$P_\mathbf{A}(\mathbf{x}) = 6x_1^2 + 7x_2^2 + 8x_3^2 + 9x_4^2 + 10x_5^2 + 11x_6^2 \\ -2x_1x_2 - 4x_1x_3 - 2x_2x_4 - 4x_3x_4 - 2x_3x_5 - 6x_4x_6 - 10x_5x_6$$

After the splitting of the weighted graph of $\mathbf{A}$, the quadratic form of $\mathbf{A}$ is also split:

(A1.2) $$P(\tilde{\mathbf{x}}) = P(\tilde{\mathbf{x}}_1) + P(\tilde{\mathbf{x}}_2)$$

$$P(\tilde{\mathbf{x}}_1) = 6x_1^2 + 7x_2^2 + 4.8x_{3a}^2 + 3.5x_{4a}^2 - 2x_1x_2 - 4x_1x_{3a} - 2x_2x_{4a} - 1.8x_{3a}x_{4a}$$

$$P(\tilde{\mathbf{x}}_2) = 3.2x_{3b}^2 + 5.5x_{4b}^2 + 10x_5^2 + 11x_6^2 - 2.2x_{3b}x_{4b} - 2x_{3b}x_5 - 6x_{4b}x_6 - 10x_5x_6$$

Here $x_3$ is split into $x_{3a}$ and $x_{3b}$, and so is $x_4$. $P(\tilde{\mathbf{x}}_1)$ is the quadratic form of $\tilde{\mathbf{A}}_1$, and $P(\tilde{\mathbf{x}}_2)$ is the quadratic form of $\tilde{\mathbf{A}}_2$, as given in Example 4.1.

If we merge $x_{3a}$ and $x_{3b}$ back to $x_3$, and merge $x_{4a}$ and $x_{4b}$ to $x_4$,

(A1.3) $$\begin{cases} x_{3a} = x_{3b} = x_3 \\ x_{4a} = x_{4b} = x_4 \end{cases}$$

then (A1.2) is changed back to (A1.1). This indicates that the variable splitting technique is also reversible.

After introducing the conception of the variable splitting technique, we begin to prove Lemma A1.1.

First, we consider a trivial case that all the vertices are on the boundary, ie. $G_B = G_a$. Then $\mathbf{A}$ could be split into $\lambda\mathbf{A}$ and $(1-\lambda)\mathbf{A}$, where $\lambda \in (0,1)$. If $\mathbf{A}$ is SPD, then both $\lambda\mathbf{A}$ and $(1-\lambda)\mathbf{A}$ must be SPD. As the result, $\mathbf{A}$ is split into 2 SPD subgraphs.

Then, we make use of the induction method. Assume $n$ is the dimension of $\mathbf{A}$.

**Step 1.** When $n=1$, there is only one vertex in $G_e$, and this vertex must be on the boundary. This is the trivial case, so Lemma A1.1 is true when $n=1$.

**Step 2.** When $n=2$, $P_A(\mathbf{x}) = a_{11}x_1^2 + 2a_{12}x_1x_2 + a_{22}x_2^2$. To split $\mathbf{A}$ into two subgraphs, there are three and only three ways to choose the boundary vertices.

(1). If $x_2$ is the boundary, then using the method of completing the square, we get:

$$P_A(\mathbf{x}) = a_{11}(x_1 + \frac{a_{12}}{a_{11}}x_2)^2 + (\frac{a_{11}a_{22} - a_{12}^2}{a_{11}})x_2^2.$$

Since $\mathbf{A}$ is SPD, $P_A(\mathbf{x}) > 0$ holds for any $\mathbf{x}$, $\mathbf{x} \neq \mathbf{0}$, thus we have:

$a_{11} > 0$ and $\Delta = a_{11}a_{22} - a_{12}^2 > 0$.

Splitting $x_2$ into $x_{2a}$ and $x_{2b}$, we get:

$$P(\tilde{\mathbf{x}}) = \left(a_{11}(x_1 + \frac{a_{12}}{a_{11}}x_{2a})^2 + \lambda(\frac{a_{11}a_{22} - a_{12}^2}{a_{11}})x_{2a}^2\right) + \left((1-\lambda)\frac{a_{11}a_{22} - a_{12}^2}{a_{11}}x_{2b}^2\right)$$

$$= P_\alpha(x_1, x_{2a}) + P_\beta(x_1, x_{2b}), \quad \lambda \in (0,1)$$

It's easy to know that both $P_\alpha$ and $P_\beta$ are positive-definite.

The corresponding matrix of $P_\alpha$ is:

$$\tilde{\mathbf{A}}_\alpha = \begin{pmatrix} a_{11} & a_{12} \\ a_{21} & \lambda a_{22} + (1-\lambda)\frac{a_{12}^2}{a_{11}} \end{pmatrix}$$

The corresponding matrix of $P_\beta$ is a $1 \times 1$ matrix shown below:

$$\tilde{\mathbf{A}}_\beta = \left((1-\lambda)a_{22} - (1-\lambda)\frac{a_{12}^2}{a_{11}}\right)$$

So, $\mathbf{A}$ is split into $\tilde{\mathbf{A}}_\alpha$ and $\tilde{\mathbf{A}}_\beta$, both of which are SPD.

(2). If $x_1$ is the boundary, $\mathbf{A}$ is also able to be split into two SPD subgraphs, because we may swap $x_1$ and $x_2$ and the conclusion for $x_2$ is also valid for $x_1$.

(3). If both $x_1$ and $x_2$ are on the boundary, this is the trivial case which has been settled before.

As the result, we conclude that Lemma A1.1 is true when $n = 2$.

**Step 3.** Assume that Lemma A1.1 is true when $n = k - 1$.

**Step 4.** When $n = k$, we assume that there is at least one vertex which is not on the boundary; otherwise, if all the vertices are on the boundary, this is the trivial case settled before. Without loss of generality, suppose that $x_k$ is not on the boundary.

$$P_A(\mathbf{x}_k) = \mathbf{x}_k^{\mathrm{T}} \mathbf{A}_k \mathbf{x}_k = \begin{pmatrix} x_1 \\ x_2 \\ \vdots \\ x_{k-1} \\ x_k \end{pmatrix}^{\mathrm{T}} \begin{pmatrix} a_{11} & a_{12} & \cdots & a_{1(k-1)} & a_{1k} \\ a_{21} & a_{22} & \cdots & a_{2(k-1)} & a_{2k} \\ \vdots & \vdots & \ddots & \vdots & \vdots \\ a_{(k-1)1} & a_{(k-1)2} & \cdots & a_{(k-1)(k-1)} & a_{(k-1)k} \\ a_{k1} & a_{k2} & \cdots & a_{k(k-1)} & a_{kk} \end{pmatrix} \begin{pmatrix} x_1 \\ x_2 \\ \vdots \\ x_{k-1} \\ x_k \end{pmatrix}$$

$$= \begin{pmatrix} x_1 \\ x_2 \\ \vdots \\ x_{k-1} \end{pmatrix}^{\mathrm{T}} \begin{pmatrix} a_{11} & a_{12} & \cdots & a_{1(k-1)} \\ a_{21} & a_{22} & \cdots & a_{2(k-1)} \\ \vdots & \vdots & \ddots & \vdots \\ a_{(k-1)1} & a_{(k-1)2} & \cdots & a_{(k-1)(k-1)} \end{pmatrix} \begin{pmatrix} x_1 \\ x_2 \\ \vdots \\ x_{k-1} \end{pmatrix} + \sum_{i=1}^{k-1} 2 a_{ik} x_i x_k + a_{kk} x_k^2$$

$$= \mathbf{x}_{k-1}^{\mathrm{T}} \mathbf{A}_{k-1} \mathbf{x}_{k-1} + \sum_{i=1}^{k-1} 2 a_{ik} x_i x_k + a_{kk} x_k^2$$

$$= P(\mathbf{x}_{k-1}) + \sum_{i=1}^{k-1} 2 a_{ik} x_i x_k + a_{kk} x_k^2$$

We set $y_k = x_k + \sum_{i=1}^{k-1} \dfrac{a_{ik}}{a_{kk}} x_i$, so that $x_k = y_k - \sum_{i=1}^{k-1} \dfrac{a_{ik}}{a_{kk}} x_i$.

$$P_A(\mathbf{x}_k) = P(\mathbf{x}_{k-1}) + \sum_{i=1}^{k-1} 2 a_{ik} x_i (y_k - \sum_{i=1}^{k-1} \dfrac{a_{ik}}{a_{kk}} x_i) + a_{kk}(y_k - \sum_{i=1}^{k-1} \dfrac{a_{ik}}{a_{kk}} x_i)^2$$

$$= P(\mathbf{x}_{k-1}) + \sum_{i=1}^{k-1} 2 a_{ik} x_i y_k - \sum_{i=1}^{k-1} \sum_{j=1}^{k-1} 2 \dfrac{a_{ik} a_{jk}}{a_{kk}} x_i x_j + a_{kk} y_k^2 - \sum_{i=1}^{k-1} 2 a_{ik} x_i y_k + \sum_{i=1}^{k-1} \sum_{j=1}^{k-1} \dfrac{a_{ik} a_{jk}}{a_{kk}} x_i x_j$$

$$= P(\mathbf{x}_{k-1}) - \sum_{i=1}^{k-1} \sum_{j=1}^{k-1} \dfrac{a_{ik} a_{jk}}{a_{kk}} x_i x_j + a_{kk} y_k^2$$

$$= \hat{P}(\mathbf{x}_{k-1}) + a_{kk} y_k^2$$

where $\hat{P}(\mathbf{x}_{k-1}) = P(\mathbf{x}_{k-1}) - \sum_{i=1}^{k-1} \sum_{j=1}^{k-1} \dfrac{a_{ik} a_{jk}}{a_{kk}} x_i x_j$.

Because $P_A(\mathbf{x}_k)$ is positive-definite, $\forall \mathbf{x}_k \in \mathbb{R}^k$, $\mathbf{x}_k \neq \mathbf{0}$, we set $x_k = -\sum_{i=1}^{k-1} \dfrac{a_{ik}}{a_{kk}} x_i$, $\mathbf{x}_k = (\mathbf{x}_{k-1}, x_k)^{\mathrm{T}}$, then we have $y_k = 0$ and $\hat{P}(\mathbf{x}_{k-1}) = P_A(\mathbf{x}_k) > 0$, $\forall \mathbf{x}_{k-1} \in \mathbb{R}^{k-1}$, $\mathbf{x}_{k-1} \neq \mathbf{0}$. Therefore $\hat{P}(\mathbf{x}_{k-1})$ is positive-definite. As $\hat{P}(\mathbf{x}_{k-1})$ is a quadratic form of ($k-1$) dimensions, it could be arbitrarily split into two positive-definite quadratic forms, as assumed in Step 3.

$$\hat{P}(\tilde{\mathbf{x}}_{k-1}) = P(\tilde{\mathbf{x}}_\alpha) + P(\tilde{\mathbf{x}}_\beta)$$

This means that the electric graph of $\hat{P}(\mathbf{x}_{k-1})$ is split into two SPD subgraphs, $\tilde{G}_\alpha$ and $\tilde{G}_\beta$.

We know that $x_k$ should not be connected to both $\tilde{G}_\alpha$ and $\tilde{G}_\beta$, because $x_k$ is not on the boundary, as we assumed at the beginning of Step 4. Without loss of generality, assume that $x_k$ is connected to $\tilde{G}_\alpha$. Then,

$$P_A(\tilde{\mathbf{x}}_k) = \hat{P}(\tilde{\mathbf{x}}_{k-1}) + a_{kk} y_k^2$$

$$= P(\tilde{\mathbf{x}}_\alpha) + P(\tilde{\mathbf{x}}_\beta) + a_{kk} y_k^2$$

$$= \left( P(\tilde{\mathbf{x}}_\alpha) + a_{kk} y_k^2 \right) + P(\tilde{\mathbf{x}}_\beta)$$

$$= \left( P(\tilde{\mathbf{x}}_\alpha) + a_{kk}(x_k + \sum_{i=1}^{k-1} \frac{a_{ik}}{a_{kk}} x_i)^2 \right) + P(\tilde{\mathbf{x}}_\beta)$$

$$= \left( P(\tilde{\mathbf{x}}_\alpha) + a_{kk}(x_k + \sum_{x_i \in \tilde{G}_\alpha} \frac{a_{ik}}{a_{kk}} x_i)^2 \right) + P(\tilde{\mathbf{x}}_\beta)$$

$$= P(\tilde{\mathbf{x}}_\alpha, x_k) + P(\tilde{\mathbf{x}}_\beta)$$

So, $P_A(\mathbf{x}_k)$ has been split into $P(\tilde{\mathbf{x}}_\alpha, x_k)$ and $P(\tilde{\mathbf{x}}_\beta)$ by variable partitioning, and both of them are positive-definite. According to Lemma A1.4, Lemma A1.1 is true when $n = k$.

**Step 5.** We conclude that Lemma A1.1 is true for arbitrary $n$.

As long as Lemma A1.1 is proved, it is straightforward to prove the Conformal Splitting existence theorem (Theorem 4.2), since one SPD graph could be split for ($N-1$) times to get $N$ SPD subgraphs.

∎

## Appendix 2. Proof for the reversibility theorem.

In Section 4 we have introduced the Electric Vertex Splitting technique from the viewpoint of a local subgragh; however, this local viewpoint is not suited to prove the reversibility theory. What we need is a global viewpoint for this splitting technique, which is presented here. The relationship between the global viewpoint and the local viewpoint is also discussed. And then, we give a basic proof for the reversibility theory (Theorem 4.1). All the discussion is bounded to the level-one splitting technique.

In Section 4, the electric graph $G_e$ of $\mathbf{Ax} = \mathbf{b}$ has been partitioned into $N$ separated subgraghs, $M_j, j = 1, 2, \cdots, N$, and each subgragh could be described by (4.4). Then, we define:

$$\tilde{\mathbf{A}} = \begin{bmatrix} \tilde{\mathbf{A}}_1 & & & 0 \\ & \tilde{\mathbf{A}}_2 & & \\ & & \ddots & \\ 0 & & & \tilde{\mathbf{A}}_N \end{bmatrix}, \quad \tilde{\mathbf{x}} = \begin{bmatrix} \mathbf{x}_1 \\ \mathbf{x}_2 \\ \vdots \\ \mathbf{x}_N \end{bmatrix}, \quad \tilde{\mathbf{b}} = \begin{bmatrix} \tilde{\mathbf{b}}_1 \\ \tilde{\mathbf{b}}_2 \\ \vdots \\ \tilde{\mathbf{b}}_N \end{bmatrix}, \quad \tilde{\boldsymbol{\omega}} = \begin{bmatrix} \tilde{\boldsymbol{\omega}}_1 \\ \tilde{\boldsymbol{\omega}}_2 \\ \vdots \\ \tilde{\boldsymbol{\omega}}_N \end{bmatrix}$$

As the result, the split system could be expressed by:

(A2.1) $$\tilde{\mathbf{A}}\tilde{\mathbf{x}} = \tilde{\mathbf{b}} + \tilde{\boldsymbol{\omega}}$$

Here $\tilde{\mathbf{A}}$ is called the split matrix of $\mathbf{A}$. (A2.1) is called the split system of the original system $\mathbf{Ax} = \mathbf{b}$. However, (A2.1) is still not suited to express the proof. We need another way to achieve this.

We define $\Gamma_{boundary}$ to be an ordered set of all the boundary vertices, and $\Gamma_{inner}$ an ordered set including all the inner vertices. Further, we define $\mathbf{u}$ the voltage vector corresponding to $\Gamma_{boundary}$, and $\mathbf{y}$ the voltage vector of $\Gamma_{inner}$.

As the result, the original linear system $\mathbf{Ax} = \mathbf{b}$ could be reformatted into (A2.2):

(A2.2) $$\begin{bmatrix} \mathbf{C} & \mathbf{E} \\ \mathbf{F} & \mathbf{D} \end{bmatrix} \begin{bmatrix} \mathbf{u} \\ \mathbf{y} \end{bmatrix} = \begin{bmatrix} \mathbf{f} \\ \mathbf{g} \end{bmatrix}$$

Then, we partition the electric graph of this system using the Electric Vertex Splitting technique, and every boundary vertex is split into a pair of twin vertices, one of which is called the senior vertex, and the other is called the junior vertex.

Hence, we define $\Gamma_{se}$ to be an ordered set of all the senior vertices, and $\Gamma_{ju}$ an ordered set of all the junior vertices. The orders of $\Gamma_{se}$, $\Gamma_{ju}$ and $\Gamma_{boundary}$ are accordant. Then, we define $\mathbf{u}_{se}$ to be the corresponding voltage vector of $\Gamma_{se}$, and $\mathbf{u}_{ju}$ the voltage vector of $\Gamma_{ju}$.

Consequently, (A2.2) is split to (A2.3).

$$\text{(A2.3)} \quad \begin{bmatrix} \mathbf{C}_{se} & \mathbf{0} & \mathbf{E}_{se} \\ \mathbf{0} & \mathbf{C}_{ju} & \mathbf{E}_{ju} \\ \mathbf{F}_{se} & \mathbf{F}_{ju} & \mathbf{D} \end{bmatrix} \begin{bmatrix} \mathbf{u}_{se} \\ \mathbf{u}_{ju} \\ \mathbf{y} \end{bmatrix} = \begin{bmatrix} \mathbf{f}_{se} \\ \mathbf{f}_{ju} \\ \mathbf{g} \end{bmatrix} + \begin{bmatrix} \boldsymbol{\omega}_{se} \\ \boldsymbol{\omega}_{ju} \\ \mathbf{0} \end{bmatrix}$$

where $\mathbf{C}_{se} + \mathbf{C}_{ju} = \mathbf{C}$, $\mathbf{E}_{se} + \mathbf{E}_{ju} = \mathbf{E}$, $\mathbf{f}_{se} + \mathbf{f}_{ju} = \mathbf{f}$. These three equations give a straightforward explanation of the splitting of the vertex weights, current sources and edge weights of the boundary vertices in Section 4. $\boldsymbol{\omega}_{se}$ and $\boldsymbol{\omega}_{ju}$ are the inflow currents of $\Gamma_{ju}$ and $\Gamma_{se}$, respectively.

(A2.3) could be represented by (A2.4), for short.

$$\text{(A2.4)} \quad \overline{\mathbf{A}}\overline{\mathbf{x}} = \overline{\mathbf{b}} + \overline{\boldsymbol{\omega}}$$

Here $\overline{\mathbf{A}}$ is symmetric.

**Lemma A2.1**: $\overline{\mathbf{A}}$ is a reordering of $\tilde{\mathbf{A}}$, and (A2.4) is equivalent to (A2.1).

Proof: The original electric graphs of $\overline{\mathbf{A}}$ and $\tilde{\mathbf{A}}$ are same, the splitting manners to generate them are same as well, and the only difference is the ordering of the unknowns in $\overline{\mathbf{x}}$ and $\tilde{\mathbf{x}}$, so Lemma A2.1 is right.

∎

Then we are going to prove the reversibility theorem, which could be re-expressed as Lemma A2.2:

**Lemma A2.2**: If $\mathbf{u}_{se} = \mathbf{u}_{ju}$, $\boldsymbol{\omega}_{se} = -\boldsymbol{\omega}_{ju}$, then $\overline{\mathbf{A}}\overline{\mathbf{x}} = \overline{\mathbf{b}} + \overline{\boldsymbol{\omega}}$ becomes $\mathbf{A}\mathbf{x} = \mathbf{b}$.

Proof: Set $\mathbf{u}_{se} = \mathbf{u}_{ju} = \mathbf{u}$, $\boldsymbol{\omega}_{se} = -\boldsymbol{\omega}_{ju} = \boldsymbol{\omega}$, then:

$$\begin{bmatrix} \mathbf{C}_{se} & 0 & \mathbf{E}_{se} \\ 0 & \mathbf{C}_{ju} & \mathbf{E}_{ju} \\ \mathbf{F}_{se} & \mathbf{F}_{ju} & \mathbf{D} \end{bmatrix} \begin{bmatrix} \mathbf{u} \\ \mathbf{u} \\ \mathbf{y} \end{bmatrix} = \begin{bmatrix} \mathbf{f}_{se} \\ \mathbf{f}_{ju} \\ \mathbf{g} \end{bmatrix} + \begin{bmatrix} \boldsymbol{\omega} \\ -\boldsymbol{\omega} \\ \mathbf{0} \end{bmatrix}$$

Eliminate $\boldsymbol{\omega}$:

$$\begin{bmatrix} \mathbf{C}_{se} + \mathbf{C}_{ju} & \mathbf{E}_{se} + \mathbf{E}_{ju} \\ \mathbf{F}_{se} + \mathbf{F}_{ju} & \mathbf{D} \end{bmatrix} \begin{bmatrix} \mathbf{u} \\ \mathbf{y} \end{bmatrix} = \begin{bmatrix} \mathbf{f}_{se} + \mathbf{f}_{ju} \\ \mathbf{g} \end{bmatrix}$$

Because $\mathbf{C}_{se} + \mathbf{C}_{ju} = \mathbf{C}$, $\mathbf{E}_{se} + \mathbf{E}_{ju} = \mathbf{E}$, $\mathbf{f}_{se} + \mathbf{f}_{ju} = \mathbf{f}$, then we get (A2.2), which is $\mathbf{Ax} = \mathbf{b}$.

∎

Finally, we present Lemma A2.3, which will be useful to prove the convergence theorem in Appendix 3.

**Lemma A2.3**: If there exists a partition scheme which assures that $\tilde{\mathbf{A}}_j$ is SPD, $j = 1, 2, \cdots N$, then $\tilde{\mathbf{A}}$ is SPD, and $\bar{\mathbf{A}}$ is SPD, consequently.

This conclusion is straightforward and the proof is omitted.

The above mathematical description of the Electric Vertex Splitting technique is only for the level-one splitting technique. The cases for the multilevel splitting techniques will be more complex and will be given in the next edition of this paper.

## Appendix 3. Proof for the convergence theorem.

Here we give a basic proof for the convergence theory (Theorem 6.1) of VTM. We only focus on the level-one splitting technique.

Assume the original graph $G_e$ is partitioned into $N$ separated subgraghs, $M_j, j = 1, 2, \cdots, N$, following some partition scheme. $G_e$ is SPD, and all the subgragh are SPD, i.e. $\tilde{\mathbf{A}}_j$ is SPD, $j = 1, 2, \cdots N$.

As described in Section 5, we add one VTLs between each pair of twin vertices. Based on the global view introduced in Appendix 2, we have,

(A3.1) $$\begin{cases} \mathbf{u}_{ju}^k + \mathbf{Z}\boldsymbol{\omega}_{ju}^k = \mathbf{u}_{se}^{k-1} - \mathbf{Z}\boldsymbol{\omega}_{se}^{k-1} \\ \mathbf{u}_{se}^k + \mathbf{Z}\boldsymbol{\omega}_{se}^k = \mathbf{u}_{ju}^{k-1} - \mathbf{Z}\boldsymbol{\omega}_{ju}^{k-1} \end{cases}$$

where **Z** should be SPD. We call **Z** the global characteristic impedance matrix of the VTLs. If all the local characteristic matrices $\mathbf{Z}_j$, $j = 1, 2, \cdots, N$, are positive diagonal matrices, then **Z** is a positive diagonal matrix as well.

Define $\tilde{\mathbf{Z}} = diag(\mathbf{Z}_1, \mathbf{Z}_2, \cdots \mathbf{Z}_N)$, and $\mathbf{M} = \begin{bmatrix} \mathbf{Z} & \mathbf{0} \\ \mathbf{0} & \mathbf{Z} \end{bmatrix}$, then, we have:

**Lemma A3.1**: **M** is a reordering of $\tilde{\mathbf{Z}}$.

Proof: **M** and $\tilde{\mathbf{Z}}$ are different ways to express the characteristic impedances of the VTLs, so they are equivalent.

∎

Remove the inner voltage **y** from (A2.3) and we get:

(A3.2) $\begin{bmatrix} \mathbf{C}_{se} - \mathbf{E}_{se}\mathbf{D}^{-1}\mathbf{F}_{se} & -\mathbf{E}_{se}\mathbf{D}^{-1}\mathbf{F}_{ju} \\ -\mathbf{E}_{ju}\mathbf{D}^{-1}\mathbf{F}_{se} & \mathbf{C}_{ju} - \mathbf{E}_{ju}\mathbf{D}^{-1}\mathbf{F}_{ju} \end{bmatrix} \begin{bmatrix} \mathbf{u}_{se} \\ \mathbf{u}_{ju} \end{bmatrix} = \begin{bmatrix} \mathbf{f}_{se} - \mathbf{E}_{se}\mathbf{D}^{-1}\mathbf{g} \\ \mathbf{f}_{ju} - \mathbf{E}_{ju}\mathbf{D}^{-1}\mathbf{g} \end{bmatrix} + \begin{bmatrix} \boldsymbol{\omega}_{se} \\ \boldsymbol{\omega}_{ju} \end{bmatrix}$

Then simplify (A3.2) into (A3.3).

(A3.3) $\mathbf{S}\hat{\mathbf{u}} = \boldsymbol{\beta} + \hat{\boldsymbol{\omega}}$

where $\mathbf{S} = \begin{bmatrix} \mathbf{C}_{se} & \\ & \mathbf{C}_{ju} \end{bmatrix} - \begin{bmatrix} \mathbf{E}_{se} \\ \mathbf{E}_{ju} \end{bmatrix} \cdot \mathbf{D}^{-1} \cdot \begin{bmatrix} \mathbf{F}_{se} & \mathbf{F}_{ju} \end{bmatrix}$.

According to Lemma A2.3, if $\tilde{\mathbf{A}}_j$ is SPD, $j = 1, 2, \cdots N$, then $\bar{\mathbf{A}}$ is SPD. Thereafter, we have Lemma A3.2.

**Lemma A3.2**: If $\bar{\mathbf{A}}$ is SPD, then **S** is SPD.

Proof: $\bar{\mathbf{A}} = \begin{bmatrix} \mathbf{C}_{se} & \mathbf{0} & \mathbf{E}_{se} \\ \mathbf{0} & \mathbf{C}_{ju} & \mathbf{E}_{ju} \\ \mathbf{F}_{se} & \mathbf{F}_{ju} & \mathbf{D} \end{bmatrix}$, as presented in Appendix 2. If $\bar{\mathbf{A}}$ is SPD, then

$\bar{\mathbf{x}}^T \bar{\mathbf{A}} \bar{\mathbf{x}} > 0$, $\forall \bar{\mathbf{x}}, \bar{\mathbf{x}} \neq \mathbf{0}$, which means that:

$$\begin{bmatrix} \mathbf{u}_{se}^T & \mathbf{u}_{ju}^T & \mathbf{y}^T \end{bmatrix} \begin{bmatrix} \mathbf{C}_{se} & \mathbf{0} & \mathbf{E}_{se} \\ \mathbf{0} & \mathbf{C}_{ju} & \mathbf{E}_{ju} \\ \mathbf{F}_{se} & \mathbf{F}_{ju} & \mathbf{D} \end{bmatrix} \begin{bmatrix} \mathbf{u}_{se} \\ \mathbf{u}_{ju} \\ \mathbf{y} \end{bmatrix} > 0,$$

$\forall \mathbf{u}_{se} \neq \mathbf{0}, \forall \mathbf{u}_{ju} \neq \mathbf{0}, \forall \mathbf{y} \neq \mathbf{0}$.

If we set $\mathbf{y} = -\mathbf{D}^{-1}\mathbf{F}_{se}\mathbf{u}_{se} - \mathbf{D}^{-1}\mathbf{F}_{ju}\mathbf{u}_{ju}$, then

(A3.4) $\begin{bmatrix} \mathbf{u}_{se}^T & \mathbf{u}_{ju}^T & -\mathbf{u}_{se}^T\mathbf{E}_{se}\mathbf{D}^{-1} - \mathbf{u}_{ju}^T\mathbf{E}_{ju}\mathbf{D}^{-1} \end{bmatrix} \begin{bmatrix} \mathbf{C}_{se} & 0 & \mathbf{E}_{se} \\ 0 & \mathbf{C}_{ju} & \mathbf{E}_{ju} \\ \mathbf{F}_{se} & \mathbf{F}_{ju} & \mathbf{D} \end{bmatrix} \begin{bmatrix} \mathbf{u}_{se} \\ \mathbf{u}_{ju} \\ -\mathbf{D}^{-1}\mathbf{F}_{se}\mathbf{u}_{se} - \mathbf{D}^{-1}\mathbf{F}_{ju}\mathbf{u}_{ju} \end{bmatrix} > 0$,

$\forall \mathbf{u}_{se} \neq \mathbf{0}$, $\forall \mathbf{u}_{ju} \neq \mathbf{0}$.

(A3.4) could be written as:

(A3.5) $\begin{bmatrix} \mathbf{u}_{se}^T & \mathbf{u}_{ju}^T \end{bmatrix} \left( \begin{bmatrix} \mathbf{C}_{se} & \\ & \mathbf{C}_{ju} \end{bmatrix} - \begin{bmatrix} \mathbf{E}_{se} \\ \mathbf{E}_{ju} \end{bmatrix} \cdot \mathbf{D}^{-1} \cdot \begin{bmatrix} \mathbf{F}_{se} & \mathbf{F}_{ju} \end{bmatrix} \right) \begin{bmatrix} \mathbf{u}_{se} \\ \mathbf{u}_{ju} \end{bmatrix} > 0$,

$\forall \mathbf{u}_{se} \neq \mathbf{0}$, $\forall \mathbf{u}_{ju} \neq \mathbf{0}$.

Because $\mathbf{S} = \begin{bmatrix} \mathbf{C}_{se} & \\ & \mathbf{C}_{ju} \end{bmatrix} - \begin{bmatrix} \mathbf{E}_{se} \\ \mathbf{E}_{ju} \end{bmatrix} \cdot \mathbf{D}^{-1} \cdot \begin{bmatrix} \mathbf{F}_{se} & \mathbf{F}_{ju} \end{bmatrix}$, (A3.5) could be expresses as:

$$\hat{\mathbf{u}}^T \mathbf{S} \hat{\mathbf{u}} > 0, \quad \forall \hat{\mathbf{u}}, \ \hat{\mathbf{u}} \neq \mathbf{0}.$$

As the result, $\mathbf{S}$ is SPD.

∎

Reformat (A3.1) into a totally matrix-vector form and we get:

$$\begin{bmatrix} \mathbf{u}_{se}^k \\ \mathbf{u}_{ju}^k \end{bmatrix} + \begin{bmatrix} \mathbf{Z} & 0 \\ 0 & \mathbf{Z} \end{bmatrix} \begin{bmatrix} \boldsymbol{\omega}_{se}^k \\ \boldsymbol{\omega}_{ju}^k \end{bmatrix} = \begin{bmatrix} \mathbf{u}_{ju}^{k-1} \\ \mathbf{u}_{se}^{k-1} \end{bmatrix} - \begin{bmatrix} \mathbf{Z} & 0 \\ 0 & \mathbf{Z} \end{bmatrix} \begin{bmatrix} \boldsymbol{\omega}_{ju}^{k-1} \\ \boldsymbol{\omega}_{se}^{k-1} \end{bmatrix}$$

Define the row exchange matrix $\mathbf{J} = \begin{bmatrix} 0 & \mathbf{I} \\ \mathbf{I} & 0 \end{bmatrix}$, where $\mathbf{I}$ is the identity matrix.

(A3.6) $\begin{bmatrix} \mathbf{u}_{se}^k \\ \mathbf{u}_{ju}^k \end{bmatrix} + \begin{bmatrix} \mathbf{Z} & 0 \\ 0 & \mathbf{Z} \end{bmatrix} \begin{bmatrix} \boldsymbol{\omega}_{se}^k \\ \boldsymbol{\omega}_{ju}^k \end{bmatrix} = \begin{bmatrix} 0 & \mathbf{I} \\ \mathbf{I} & 0 \end{bmatrix} \begin{bmatrix} \mathbf{u}_{se}^{k-1} \\ \mathbf{u}_{ju}^{k-1} \end{bmatrix} - \begin{bmatrix} 0 & \mathbf{I} \\ \mathbf{I} & 0 \end{bmatrix} \begin{bmatrix} \mathbf{Z} & 0 \\ 0 & \mathbf{Z} \end{bmatrix} \begin{bmatrix} \boldsymbol{\omega}_{se}^{k-1} \\ \boldsymbol{\omega}_{ju}^{k-1} \end{bmatrix}$

(A3.6) could be simply expressed by (A3.7).

(A3.7) $\hat{\mathbf{u}}^k + \mathbf{M}\hat{\boldsymbol{\omega}}^k = \mathbf{J}(\hat{\mathbf{u}}^{k-1} - \mathbf{M}\hat{\boldsymbol{\omega}}^{k-1})$

Remind that $\mathbf{M} = \begin{bmatrix} \mathbf{Z} & 0 \\ 0 & \mathbf{Z} \end{bmatrix}$. Because $\mathbf{Z}$ is SPD, $\mathbf{M}$ is SPD.

According to (A3.3), remove $\hat{\boldsymbol{\omega}}^k$ from (A3.7). We get:

$$\hat{\mathbf{u}}^k + \mathbf{M}(\mathbf{S}\hat{\mathbf{u}}^k - \boldsymbol{\beta}) = \mathbf{J}\hat{\mathbf{u}}^{k-1} - \mathbf{JM}(\mathbf{S}\hat{\mathbf{u}}^{k-1} - \boldsymbol{\beta})$$

$$\hat{\mathbf{u}}^k = (\mathbf{I} + \mathbf{MS})^{-1}\mathbf{J}(\mathbf{I} - \mathbf{MS})\hat{\mathbf{u}}^{k-1} + (\mathbf{I} + \mathbf{MS})^{-1}(\mathbf{JM} + \mathbf{M})\boldsymbol{\beta}$$

Let $\mathbf{P} = (\mathbf{I} + \mathbf{MS})^{-1}\mathbf{J}(\mathbf{I} - \mathbf{MS})$, $\boldsymbol{\gamma} = (\mathbf{I} + \mathbf{MS})^{-1}(\mathbf{JM} + \mathbf{M})\boldsymbol{\beta}$, then,

(A3.8) $$\hat{\mathbf{u}}^k = \mathbf{P}\hat{\mathbf{u}}^{k-1} + \boldsymbol{\gamma}$$

**Lemma A3.3**: If $\mathbf{M}$ is SPD, and $\mathbf{S}$ is SPD, then $\sqrt{\mathbf{M}}\mathbf{S}\sqrt{\mathbf{M}} = \mathbf{QTQ}^{\mathrm{T}}$. Here $\mathbf{QQ}^{\mathrm{T}} = \mathbf{I}$, $\sqrt{\mathbf{M}}\sqrt{\mathbf{M}} = \mathbf{M}$, $\mathbf{T} = diag(t_1, t_2, \cdots, t_r)$, $t_i > 0$, $i = 1, 2, \cdots, r$. $r$ is the dimension of $\mathbf{S}$.

Proof: Because $\mathbf{M}$ is SPD, and $\mathbf{S}$ is SPD, $\sqrt{\mathbf{M}}\mathbf{S}\sqrt{\mathbf{M}}$ is SPD, then there exists a real orthogonal matrix $\mathbf{Q}$ such that $\sqrt{\mathbf{M}}\mathbf{S}\sqrt{\mathbf{M}} = \mathbf{QTQ}^{\mathrm{T}}$, where $\mathbf{T}$ is a positive diagonal matrix.

∎

**Lemma A3.4**: $\mathbf{MS} = \sqrt{\mathbf{M}}\mathbf{QTQ}^{\mathrm{T}}\sqrt{\mathbf{M}}^{-1}$

Proof: $\mathbf{MS} = \sqrt{\mathbf{M}}(\sqrt{\mathbf{M}}S\sqrt{\mathbf{M}})\sqrt{\mathbf{M}}^{-1} = \sqrt{\mathbf{M}}(\mathbf{QTQ}^{\mathrm{T}})\sqrt{\mathbf{M}}^{-1}$

**Lemma A3.5**: $\sqrt{\mathbf{M}}^{-1}\mathbf{J}\sqrt{\mathbf{M}} = \mathbf{J}$

Proof: $\sqrt{\mathbf{M}}^{-1}\mathbf{J}\sqrt{\mathbf{M}} = \begin{bmatrix} \sqrt{\mathbf{Z}} & 0 \\ 0 & \sqrt{\mathbf{Z}} \end{bmatrix}^{-1} \times \begin{bmatrix} 0 & \mathbf{I} \\ \mathbf{I} & 0 \end{bmatrix} \times \begin{bmatrix} \sqrt{\mathbf{Z}} & 0 \\ 0 & \sqrt{\mathbf{Z}} \end{bmatrix}$

$$= \begin{bmatrix} \sqrt{\mathbf{Z}}^{-1} & 0 \\ 0 & \sqrt{\mathbf{Z}}^{-1} \end{bmatrix} \times \begin{bmatrix} 0 & \mathbf{I} \\ \mathbf{I} & 0 \end{bmatrix} \times \begin{bmatrix} \sqrt{\mathbf{Z}} & 0 \\ 0 & \sqrt{\mathbf{Z}} \end{bmatrix}$$

$$= \begin{bmatrix} \sqrt{\mathbf{Z}}^{-1} & 0 \\ 0 & \sqrt{\mathbf{Z}}^{-1} \end{bmatrix} \times \begin{bmatrix} 0 & \sqrt{\mathbf{Z}} \\ \sqrt{\mathbf{Z}} & 0 \end{bmatrix} = \begin{bmatrix} 0 & \mathbf{I} \\ \mathbf{I} & 0 \end{bmatrix} = \mathbf{J}$$

∎

According to Lemma A3.4 and A3.5, we write,

$$\mathbf{P} = \left(\mathbf{I} + \sqrt{\mathbf{M}}\mathbf{Q}\mathbf{T}\mathbf{Q}^T\sqrt{\mathbf{M}}^{-1}\right)^{-1}\mathbf{J}\left(\mathbf{I} - \sqrt{\mathbf{M}}\mathbf{Q}\mathbf{T}\mathbf{Q}^T\sqrt{\mathbf{M}}^{-1}\right)$$

$$= \left(\sqrt{\mathbf{M}}\mathbf{Q}\mathbf{Q}^T\sqrt{\mathbf{M}}^{-1} + \sqrt{\mathbf{M}}\mathbf{Q}\mathbf{T}\mathbf{Q}^T\sqrt{\mathbf{M}}^{-1}\right)^{-1}\mathbf{J}\left(\sqrt{\mathbf{M}}\mathbf{Q}\mathbf{Q}^T\sqrt{\mathbf{M}}^{-1} - \sqrt{\mathbf{M}}\mathbf{Q}\mathbf{T}\mathbf{Q}^T\sqrt{\mathbf{M}}^{-1}\right)$$

$$= \sqrt{\mathbf{M}}\mathbf{Q}(\mathbf{I}+\mathbf{T})^{-1}\mathbf{Q}^T\sqrt{\mathbf{M}}^{-1}\mathbf{J}\sqrt{\mathbf{M}}\mathbf{Q}(\mathbf{I}-\mathbf{T})\mathbf{Q}^T\sqrt{\mathbf{M}}^{-1}$$

$$= \sqrt{\mathbf{M}}\mathbf{Q}(\mathbf{I}+\mathbf{T})^{-1}\mathbf{Q}^T\mathbf{J}\mathbf{Q}(\mathbf{I}-\mathbf{T})\mathbf{Q}^T\sqrt{\mathbf{M}}^{-1}$$

$$\mathbf{P}^k = \left(\sqrt{\mathbf{M}}\mathbf{Q}(\mathbf{I}+\mathbf{T})^{-1}\mathbf{Q}^T\mathbf{J}\mathbf{Q}(\mathbf{I}-\mathbf{T})\mathbf{Q}^T\sqrt{\mathbf{M}}^{-1}\right)^k$$

$$= \sqrt{\mathbf{M}}\mathbf{Q}(\mathbf{I}+\mathbf{T})^{-1}\mathbf{Q}^T\left(\mathbf{J}\mathbf{Q}(\mathbf{I}-\mathbf{T})(\mathbf{I}+\mathbf{T})^{-1}\mathbf{Q}^T\right)^{k-1}\mathbf{J}\mathbf{Q}(\mathbf{I}-\mathbf{T})\mathbf{Q}^T\sqrt{\mathbf{M}}^{-1}$$

Therefore,

$$\left\|\mathbf{P}^k\right\|_2 = \left\|\sqrt{\mathbf{M}}\mathbf{Q}(\mathbf{I}+\mathbf{T})^{-1}\mathbf{Q}^T\left(\mathbf{J}\mathbf{Q}(\mathbf{I}-\mathbf{T})(\mathbf{I}+\mathbf{T})^{-1}\mathbf{Q}^T\right)^{k-1}\mathbf{J}\mathbf{Q}(\mathbf{I}-\mathbf{T})\mathbf{Q}^T\sqrt{\mathbf{M}}^{-1}\right\|_2$$

$$\leq \left\|\sqrt{\mathbf{M}}\right\|\times\|\mathbf{Q}\|\times\left\|(\mathbf{I}+\mathbf{T})^{-1}\right\|\times\|\mathbf{Q}^T\|\times\left\|\left(\mathbf{J}\mathbf{Q}(\mathbf{I}-\mathbf{T})(\mathbf{I}+\mathbf{T})^{-1}\mathbf{Q}^T\right)^{k-1}\right\|$$

$$\times\|\mathbf{J}\|\times\|\mathbf{Q}\|\times\|\mathbf{I}-\mathbf{T}\|\times\|\mathbf{Q}^T\|\times\left\|\sqrt{\mathbf{M}}^{-1}\right\|$$

$$= \left\|\sqrt{\mathbf{M}}\right\|\times\left\|(\mathbf{I}+\mathbf{T})^{-1}\right\|\times\left(\|\mathbf{J}\|\times\|\mathbf{Q}\|\times\|(\mathbf{I}-\mathbf{T})(\mathbf{I}+\mathbf{T})^{-1}\|\times\|\mathbf{Q}^T\|\right)^{k-1}\times\|\mathbf{I}-\mathbf{T}\|\times\left\|\sqrt{\mathbf{M}}^{-1}\right\|$$

$$= \left\|\sqrt{\mathbf{M}}\right\|\times\left\|(\mathbf{I}+\mathbf{T})^{-1}\right\|\times\left\|(\mathbf{I}-\mathbf{T})(\mathbf{I}+\mathbf{T})^{-1}\right\|^{k-1}\times\|\mathbf{I}-\mathbf{T}\|\times\left\|\sqrt{\mathbf{M}}^{-1}\right\|$$

$$= \left\|\sqrt{\mathbf{M}}\right\|\times\left\|\sqrt{\mathbf{M}}^{-1}\right\|\times\left\|diag(\frac{1-t_1}{1+t_1},\frac{1-t_2}{1+t_2},\cdots,\frac{1-t_r}{1+t_r})\right\|^{k-1}\times\|\mathbf{I}-\mathbf{T}\|\times\left\|(\mathbf{I}+\mathbf{T})^{-1}\right\|$$

$$\leq \left\|\sqrt{\mathbf{M}}\right\|\times\left\|\sqrt{\mathbf{M}}^{-1}\right\|\times\|\mathbf{I}-\mathbf{T}\|\times\left\|(\mathbf{I}+\mathbf{T})^{-1}\right\|\times\left\{\max\left(\frac{1-t_1}{1+t_1},\frac{1-t_2}{1+t_2},\cdots,\frac{1-t_r}{1+t_r}\right)\right\}^{k-1}$$

When $k$ is large enough, $\left\|\mathbf{P}^k\right\|_2 < 1$. Then the iteration (A3.8) converges for any initial $\hat{\mathbf{u}}^0$. So we conclude that VTM is convergent.

Finally, we are going to prove that the converging result is the answer of the original system $\mathbf{Ax} = \mathbf{b}$.

Suppose that $\lim_{k \to \infty} \hat{\mathbf{u}}^k = \hat{\mathbf{U}}$, $\lim_{k \to \infty} \hat{\boldsymbol{\omega}}^k = \hat{\boldsymbol{\Omega}}$, which is equal to,

$$\begin{bmatrix} \mathbf{u}_{se}^k \\ \mathbf{u}_{ju}^k \end{bmatrix} \xrightarrow{k \to \infty} \begin{bmatrix} \mathbf{U}_{se} \\ \mathbf{U}_{ju} \end{bmatrix}, \quad \begin{bmatrix} \boldsymbol{\omega}_{se}^k \\ \boldsymbol{\omega}_{ju}^k \end{bmatrix} \xrightarrow{k \to \infty} \begin{bmatrix} \boldsymbol{\Omega}_{se} \\ \boldsymbol{\Omega}_{ju} \end{bmatrix}.$$

Then we get (A3.11) from (A3.7):

(A3.11) $\qquad \hat{\mathbf{U}} + \mathbf{M}\hat{\boldsymbol{\Omega}} = \mathbf{J}\hat{\mathbf{U}} - \mathbf{J}\mathbf{M}\hat{\boldsymbol{\Omega}}$

Multiply both sides of (A3.11) by $\mathbf{J}$. We obtain,

$$\mathbf{J}\hat{\mathbf{U}} + \mathbf{J}\mathbf{M}\hat{\boldsymbol{\Omega}} = \mathbf{J}\mathbf{J}\hat{\mathbf{U}} - \mathbf{J}\mathbf{J}\mathbf{M}\hat{\boldsymbol{\Omega}} = \hat{\mathbf{U}} - \mathbf{M}\hat{\boldsymbol{\Omega}}$$

(A3.12) $\qquad \hat{\mathbf{U}} - \mathbf{M}\hat{\boldsymbol{\Omega}} = \mathbf{J}\hat{\mathbf{U}} + \mathbf{J}\mathbf{M}\hat{\boldsymbol{\Omega}}$

Add (A3.11) to (A3.12), we get,

$$\mathbf{U} = \mathbf{J} \times \mathbf{U}$$

Thus,

$$\mathbf{I} = -\mathbf{J} \times \mathbf{I}$$

As the result,

$$\mathbf{U}_{se} = \mathbf{U}_{ju}, \quad \boldsymbol{\Omega}_{se} = -\boldsymbol{\Omega}_{ju}.$$

According to the reversibility theorem (Theorem 4.1), we conclude that the convergent result is exactly the answer to the original system. So we have proved the convergence theorem of VTM.

It should be noted that the above proof does not cover the case when there exists multilevel splitting during graph partitioning. A full proof will be given in the next edition of this paper.